\documentclass[reqno]{amsart}
\usepackage{amsmath}
\usepackage{amstext}
\usepackage{amsfonts}
\usepackage[mathscr]{eucal}
\usepackage{amscd}
\usepackage{latexsym}
\usepackage{amssymb}
\usepackage{diagrams}

\newtheorem{theorem}[]{Theorem}
\newtheorem{lemma}[theorem]{Lemma}
\newtheorem{proposition}[theorem]{Proposition}
\newtheorem{definition}[theorem]{Definition}
\newtheorem{corollary}[theorem]{Corollary}
\newtheorem{remark}[theorem]{Remark}
\theoremstyle{definition}

\newcommand{\zbb}{\mathbb{Z}}
\newcommand{\gbb}{\mathbb{G}}
\newcommand{\bsy}[1]{\boldsymbol{#1}}
\newcommand{\mi}{\text{-}}
\newcommand{\ext}{\mathrm{Ext}}
\newcommand{\tor}{\mathrm{Tor}}

\newcommand{\homz}{{\mathrm{Hom}}_{\mathbb{Z}}}
\newcommand{\otmz}{{\otimes}_{\mathbb{Z}}}
\newcommand{\enpr}{\hfill $\Box $}
\newcommand{\supar}[1]{\overset{#1}\rightarrow}
\newcommand{\suparle}[1]{\overset{#1}\leftarrow}

\newcommand{\prf}{\noindent\textbf{Proof. }}

\newcommand{\rw}{\rightarrow}

\newcommand{\xrw}{\xrightarrow}

\newcommand{\id}{\mathrm{id}}
\newcommand{\der}{\mathrm{Der}}
\newcommand{\hm}{\mathrm{Hom}}
\newcommand{\md}[1]{${#1}$-module}
\newcommand{\ifr}{\mathfrak{I}}

\newcommand{\diff}{\mathrm{Diff}\,}
\newcommand{\hund}[1]{{\hm}_{#1}}
\newcommand{\otund}[1]{{\otimes}_{#1}}

\newcommand{\pr}{\mathrm{pr}}
\newcommand{\ab}{\mathbf{Ab}}
\newcommand{\cm}{\mathbf{CM}}
\newcommand{\tgmu}{(T,G,\mu)}

\newcommand{\ccl}{\mathcal{C}}
\newcommand{\gcl}{\mathcal{G}}
\newcommand{\remb}{\rightarrowtail}
\newcommand{\ronto}{\twoheadrightarrow}
\newcommand{\p}{\pi_1}
\newcommand{\pp}{\pi'_1}

\newcommand{\rt}{\rtimes}
\newcommand{\tot}{\mathrm{Tot\,}}
\newcommand{\set}{\mathbf{Set}}
\newcommand{\s}{\scriptstyle}
\newcommand{\D}{\displaystyle}
\newcommand{\cir}{\s\,\circ\,\D}
\newcommand{\pcm}{\bsy{\mathrm{PCM}}}
\newcommand{\biset}{\text{Simpl}^2\text{Set}}
\newarrow{Eq}{}{=}{=}{=}{}
\newarrow{LR}{<}{-}{-}{-}{>}
\newarrow{Dts}{}{.}{.}{.}{}
\newarrow{Embedline}{>}{-}{-}{-}{-}
\newarrow{Intoar}{C}{-}{-}{-}{>}
\begin{document}%

\makeatletter
\renewcommand{\descriptionlabel}[1]{\hspace\labelsep \upshape\bfseries #1}
\makeatother

\title{(Co)homology of crossed modules\\ with coefficients in a ${\p}$-module}
\author{ Simona Paoli }
\email{simona@maths.warwick.ac.uk}
\address{Mathematics Institute\\ University of Warwick\\ Coventry CV4 7AL\\
 England}



\begin{abstract}
We define a cotriple (co)homology of crossed modules with coefficients in a $\p$-module.
We prove its general properties, including the connection with the existing cotriple
theories on crossed modules. We establish the relationship with the (co)homology of the
classifying space of a crossed module and with the cohomology of groups with operators.
An example and an application are given.
\end{abstract}

%\received{...}   % receive date (for example: 11 October 1999)
%\revised{...}    % receive date
%\published{...}  % publish date
%\submitted{...}  % Name of Journal's Editor, who submitted Article

%\volumeyear{...} % Volume Year
%\volumenumber{...}  % Volume Number
%\issuenumber{...}   % Issue Number

%\startpage{1}     % Page Number of first page

\maketitle

%%%%%% New section %%%%%%%%%%%%%%%%%%%%%%%%%%%%%%%%%%%%%%%%%%%%%%%%%%%%%%%%%%%%%%
\section*{Introduction}\label{intro}
In this paper we present some developments in the (co)homology of crossed modules. In the
work of Carrasco, Cegarra and Grandjean \cite{ccg} the authors proved that the category
of crossed modules is tripleable over the category of sets, hence it is an algebraic
category; then they used the resulting cotriple to construct a (co)homology theory of
crossed modules in the spirit of the Barr and Beck theory \cite{bab1}.

Later Grandjean, Ladra and Pirashvili \cite{glp} have proved that there is an exact
homology sequence
\begin{equation}\label{intro.eq1}
  \cdots\rw H_{n+1}B\tgmu\rw\zeta H^{CCG}_n\tgmu\rw H_nG\rw H_n B\tgmu\rw\cdots
\end{equation}
which relates the integral homology of the classifying space of a crossed module and the
cotriple homology of \cite{ccg}.

The (co)homology theory of \cite{ccg} has trivial coefficients. In any algebraic category
the passage from trivial coefficients for the (co)homology theory to global or local ones
is achieved by a well known procedure which consists of taking abelian group objects in
the slice category; the cohomology of a crossed module $\Phi$ with a system of global or
local coefficients is equivalent to the cohomology in the slice category of $\id_\Phi$
with trivial coefficients. Although there is no theoretic difficulty in realizing this
passage, to achieve it in practice in a concrete algebraic context like the one of
crossed modules is not entirely trivial.

One of the first questions to consider is whether it is possible identify a manageable
class of coefficients giving rise to a (co)homology theory of crossed modules which has
interesting properties and leads to applications. This was one of our motivating
questions which we have tried to answer in this paper by concentrating our attention on
two special cases of a system of local coefficients associated to a $\p$-module, where
$\p$ is the first homotopy group of the crossed module.

Our motivation for considering $\p$-module coefficients for the cotriple (co)homolo-gy
comes from the long exact sequence (\ref{intro.eq1}). Since the (co)homology of the
classifying space of a crossed module is defined in general with $\p$-module coefficients
(see \cite{elli}) it is natural to expect that an appropriate cotriple (co)homology of
crossed modules with $\p$-module coefficients would lead to long exact (co)homology
sequences generalizing (\ref{intro.eq1}) in the homology case. One of our main results,
Theorem \ref{class.the1b}, establishes precisely this.

Given a crossed module $\Phi$ with first homotopy group $\p$, the canonical projection
$\Phi\ronto(1,\p,i)\equiv\p(\Phi)$ induces a functor
\begin{equation}\label{intro.eq2}
  (\cm/\p(\Phi))_{ab}\rw(\cm/\Phi)_{ab}.
\end{equation}
A system of local coefficients for the cohomology of a crossed module $\Phi$ could be
defined as an object of $(\cm/\p(\Phi))_{ab}$. The functor (\ref{intro.eq2}) takes this
system of local coefficients into a system of global coefficients used to compute the
(co)homology.

In this paper we work with two special cases of the above system of local coefficients;
these correspond to the split extensions of crossed modules
\begin{equation}\label{intro.eq3}
  (A,1,0)\remb(A,\p,0)\ronto(1,\p,i)
\end{equation}
and
\begin{equation}\label{intro.eq4}
  (1,A,i)\remb(1,A\rt\p,i)\ronto(1,\p,i).
\end{equation}
In Sections 2 to 6 we consider the coefficients corresponding to the split extension
(\ref{intro.eq3}), and study the corresponding (co)homology theory. After recalling some
background in Section 1, in Section 2 we introduce the (co)homology. For this purpose,
the notions of module and derivation in the sense of categories of interest \cite{orz}
are used to define a derivation functor, and a dual $\diff$ functor, from crossed modules
to abelian groups. We point out that the notions of action, extensions and semidirect
product of crossed modules were also worked out by Norrie \cite{nor}, internally in the
category of crossed modules. This has been used by Vieites and Casas \cite{rodr} to give
a different approach to derivations of crossed modules.

In Section 3 we study the case of aspherical crossed modules; these form a subcategory
isomorphic to the category of surjective group homomorphisms. We prove that our
(co)homology for the aspherical crossed module corresponding to the surjective group
homomorphism $f:G\rw G'$ is isomorphic, up to a dimension shift, to the relative group
(co)homology of the pair $(G',G)$ defined by Loday \cite{loda}. Other general properties
of the (co)homology are proved in Section 4. In Section 5 we establish the relationship
between our theory and the (co)homology of the classifying space, recovering the result
of \cite{glp} in the case of homology with integral coefficients. This result is
illustrated with an example in Section 6 where we obtain some information about the
(co)homology of the crossed module corresponding to a \md{\zbb G} $M$.

The coefficients corresponding to the split extension (\ref{intro.eq4}) are treated in
Section 7. We prove that the corresponding cotriple (co)homology coincides with the
(co)homology of the classifying space of the crossed module, up to a dimension shift of
1, and in dimensions $n>0$. An application to the cohomology of the classifying space
follows.

In the last section we elucidate the relationship between the cohomology of crossed
modules with $\p$-module coefficients introduced in Section 2 and the cohomology of
groups with operators studied in \cite{cega}. In order to study the relationship between
the two theories we establish the preliminary result, which may be of independent
interest, that when $\tgmu$ is a precrossed module the cohomology $H^*_G(T,A)$ of
\cite{cega} can be described as cohomology of precrossed modules with a system of local
coefficients.

%%%%%% New section %%%%%%%%%%%%%%%%%%%%%%%%%%%%%%%%%%%%%%%%%%%%%%%%%%%%%%%%%%%%%%%%%%%
\section*{Acknowledgments}
This paper is based on part of my thesis \cite{pao}. I am very grateful to my supervisor
Alan Robinson for his advice and encouragement and to the Engineering and Physical
Sciences Research Council for financial support. I thank the Georgian Academy of Sciences
for their hospitality while part of this work was carried out; in particular I am
grateful to Hvedri Inassaridze for suggesting the initial direction of my work  and  to
Tamar Datuashvili and Teimuraz Pirashvili for useful conversations. I would also like to
thank the referee for some helpful comments and suggestions.

%%%%%%% New section %%%%%%%%%%%%%%%%%%%%%%%%%%%%%%%%%%%%%%%%%%%%%%%%%%%%%%%%%%%%%%%%
\section{Preliminaries}\label{prel}
\subsection{Crossed modules.}\label{crossm}
Recall that a \emph{crossed module} $\Phi=\tgmu$ consists of a group homomorphism
$\mu:T\rw G$ and of an action of $G$ on $T$ such that
\begin{align*}
  \mu(\,^{g}t)=g\mu(t)g^{-1}, \qquad ^{\mu(t)}t'=tt't^{-1}
\end{align*}
for each $t,t'\in T,\;g\in G$. A \emph{homomorphism of crossed modules}
$(f_T,f_G):\tgmu\rw(T',G',\mu')$ is a pair of group homomorphisms $f_T:T\rw
T',\;\;f_G:G\rw G'$ such that $\mu'f_T=f_G\mu$ and $f_T(\,^{g}t)=\,^{f_G(g)}f_T(t)$ for
all $g\in G,\;t\in T$. We denote by $\cm$ the category of crossed modules. This category
has several equivalent descriptions.

Recall that a cat$^{1}$-group consists of a group $G$ with two endomorphisms
$d_0,d_1:G\rw G$ such that
\begin{equation}\label{prel.eq1}
d_1d_0=d_0,\qquad d_0d_1=d_1,\qquad [\ker d_0,\ker d_1]=1.
\end{equation}
A \emph{morphism of cat$^{1}$-groups} $\;(G,d_0,d_1)\,\rw\,(G',d'_0,d'_1)$ is a group
homomorphism $\;f:G\rw G'$ such that $d'_if=fd_i\;\;i=0,1$.

The category of crossed modules is equivalent to the category of cat$^{1}$-groups
\cite{lodb}. Given a crossed module $\tgmu$, the corresponding cat$^{1}$-group is
$(T\rtimes G,d_0,d_1)$, $\;d_0(t,g)=(1,g)$, $\;d_1(t,g)=(1,\mu(t)g)$ for all $(t,g)\in
T\rtimes G$.

Another description of the category of crossed modules is given by its equivalence with
the category $\mathbf{SG}_{\leq 1}$ of simplicial groups whose Moore complex has length 1
\cite{lodb}. An object of $\mathbf{SG}_{\leq 1}$ is a simplicial group $G_*$ such that
$N_i(G_*)=0$ for $i>1$ while $N_1(G_*)\neq 0$ where $N_*:\mathbf{SG}_{\leq 1}\rw\cm$ is
the Moore normalization functor. This is defined by
\begin{equation*}
  N_nG_*=\underset{i>0}{\cap}\ker(\partial_i^n)\qquad n>0
\end{equation*}
with boundary map $d:N_nG_*\rw N_{n-1}G_*$, $\;d=\partial_{{0|}_{N_nG_*}}$. There is a
functor $N_*^{-1}:\cm\rw\mathbf{SG}_{\leq 1}$ with $N_*^{-1}(N_*G)\cong G_*$ which is
given by
\begin{align*}
   & N^{-1}_{n}(T,G,\mu)=T^n\rtimes G \qquad\qquad\qquad\qquad n\geq 0, \\
   & \partial_i(t_1,\ldots ,t_n,g)=(t_1,\ldots ,\hat{t}_i,\ldots,t_n,g)\qquad 1\leq i\leq
   n, \\
   & \partial_0(t_1,\ldots ,t_n,g)=(t_2t^{-1}_1,\ldots,t_nt^{-1}_1,\mu t_1g)\\
   & s_i(t_1,\ldots ,t_n,g)=(t_1,\ldots,t_{i},0,\ldots,t_n,g).
\end{align*}

Crossed modules are algebraic models for connected spaces which have trivial homotopy
groups in dimension $n>2$, called 2-types (see for instance \cite{lodb}). To any crossed
module $(T,G,\mu)$ one can associate a connected $\mathrm{CW}$-space $B(T,G,\mu)$ called
its \emph{classifying space} with
\begin{align*}
   & \pi_1 B(T,G,\mu)\cong G/\mu(T),\quad \pi_2 B(T,G,\mu)\cong\ker\mu,\quad
   \pi_n B(T,G,\mu)\cong 0\;\; \text{for}\;\; n>2.
\end{align*}
$B\tgmu$ is defined as the classifying space of the simplicial group $N_*^{-1}\tgmu$.

The \emph{homotopy groups} of a crossed module $\tgmu$ are defined as $\p=G/\mu(T)$,
$\;\pi_2=\ker\mu$, $\;\pi_n=0$ for $n\neq 1,2$. A morphism of crossed modules is called a
\emph{weak equivalence} if it induces isomorphisms of homotopy groups.

It can be proved (see for instance \cite{lodb}) that the functor $B(\mi)$ induces an
equivalence between the homotopy category of connected 2-types and the localization of
the category of crossed modules with respect to weak equivalences.

\subsection{CCG (co)homology.} In \cite{ccg} it is proved that the category of crossed
modules is tripleable over $\mathbf{Set}$, hence it is an algebraic category. It is shown
there that the functor $\mathcal{U}:\cm\rw\mathbf{Set}$, $\;\mathcal{U}\tgmu=T\times G$
has a left adjoint $\mathcal{F}:\mathbf{Set}\rw\cm$. This is given by
\begin{equation*}
  \mathcal{F}(X)=(\overline{F(X)},F(X)*F(X),i)
\end{equation*}
where $F(X)$ is the free group on $X$, $*$ is the free product, $i$ is the inclusion,
$\overline{F(X)}$ is the kernel of the map $p_2:F(X)*F(X)\rw F(X)$ determined by
$p_2u_1=0$, $\;p_2u_2=\id$, $u_1,u_2$ being the coproduct injections.

It is proved in \cite{ccg} that $\mathcal{U}$ is tripleable. This identifies the regular
epimorphisms in $\cm$ as those homomorphisms $(f_T,f_G):\tgmu\rw(T',G',\mu')$ such that
$f_T$ and $f_G$ are surjective. Hence for each set $X$ the crossed module
$\mathcal{F}(X)$ is a projective object in $\cm$; this category has enough projectives
since any crossed module $\tgmu$ admits the projective presentation
$\mathcal{F}\mathcal{U}\tgmu\ronto\tgmu$.

Let $\gbb=\mathcal{F}\mathcal{U}$ be the cotriple arising from the pair of adjoint
functors $( \mathcal{F},\mathcal{U)}$. This cotriple is used in \cite{ccg} to define a
(co)homology theory of crossed modules as follows. Recall that given a crossed module
$\tgmu$ its abelianisation is the abelian crossed module
$\tgmu_{ab}=(T/[G,T],G/[G,G],\overline{\mu})$. For each $n\geq 1$ the \emph{$n^{th}$ CCG
homology} of $\tgmu$ is the crossed module
\begin{equation*}
  H_n^{CCG}\tgmu=H_{n-1}(\gbb_{\bullet}\tgmu_{ab}).
\end{equation*}
If $(A,B,\delta)$ is an abelian crossed module (that is $A$ and $B$ are abelian groups
and $B$ acts trivially on $A$), for each $n\geq 1$ the \emph{$n^{th}$ CCG cohomology} of
$\tgmu$ with coefficients in $(A,B,\delta)$  is the abelian group
\begin{equation*}
\begin{split}
   & H^n_{CCG}(\tgmu,(A,B,\delta))=H^{n-1}\hund{\cm}(\gbb_{\bullet}\tgmu,(A,B,\delta))\cong \\
   & \cong H^{n-1}\hund{\mathbf{Ab}\,\cm}(\gbb_{\bullet}\tgmu_{ab},(A,B,\delta)).
\end{split}
\end{equation*}

\subsection{Crossed modules as category of interest.} It is a known fact that the
category of cat$^{1}$-groups is a category of groups with operations in the sense of
\cite{po}. Recall that this consists of the following data: a category of groups with a
set of operations $\Omega=\Omega_0\cup\Omega_1\cup\Omega_2$ where $\Omega_i$ is the set
of $i$-ary operations in $\Omega$ such that the group operations of identity, inverse and
multiplication (denoted $0,-,+$) are elements of $\Omega_0,\,\Omega_1,\,\Omega_2$
respectively; one has $\Omega_0=\{0\}$ and certain compatibility conditions hold (see
\cite{po}); finally there is a set of identities $E$ which includes the group laws.

In the case of cat$^{1}$-groups, $\Omega_0=\{0\}$, $\;\Omega_1=\{-\}\cup\{d_0,d_1\}$,
$\;\Omega_2=\{+\}$ and $E$ consists of the group laws and of the identities
(\ref{prel.eq1}). The compatibility conditions in this case are that $d_0,d_1$ commute
with $+$, hence they are group homomorphisms.

In the category of cat$^{1}$-groups we therefore have the notions of singular object,
module, semidirect product, derivation. For a discussion of these notions in any category
of groups with operations see for instance \cite{dat} and \cite{po}.

A cat$^{1}$-group $(A,d_0,d_1)$ is a singular object if $A$ is an abelian group. The
corresponding crossed module is then an abelian crossed module. Given a cat$^{1}$-group
$(G,s_0,s_1)$, $\;(A,d_0,d_1)$ is a $(G,s_0,s_1)$-module if $(A,d_0,d_1)$ is singular and
there is a split extension of cat$^{1}$-groups
\begin{equation*}
  (A,d_0,d_1)\remb(Q,s'_0,s'_1)\;\,\pile{\suparle{}\\ \ronto}\;\,(G,s_0,s_1).
\end{equation*}
This split singular extension induces an action of $(G,s_0,s_1)$ on $(A,d_0,d_1)$; a
derivation $D$ from $(G,s_0,s_1)$ into $(A,d_0,d_1)$ is a group derivation from $G$ into
$A$ which commutes with the 1-ary operations $\omega\in \Omega_1\verb"\"\{-\}$, that is
such that $Ds_i=d_iD$ for $i=0,1$.

Since $\cm$ is equivalent to the category of cat$^{1}$-groups, it can be considered
itself as a category of groups with operations. Moreover, since $\cm$ is tripleable over
$\mathbf{Set}$ and the set $\Omega_2$ of 2-ary operations just consists of group
multiplication, it is in fact a category of interest in the sense of \cite{orz}. In this
paper we shall use the interpretation in terms of extensions of the first and second
cotriple cohomology in a category of interest given in \cite{vale}.

\subsection{Crossed modules in the category of cat$\bsy{^1}$-groups.}\label{catone}
Let $\ccl^1\gcl$ denote the category of cat$^1$-groups and $\ccl^2\gcl$ the category of
cat$^2$-groups; we refer to \cite{lodb} for the definition of cat$^2$-group. Since
$\ccl^1\gcl$ is a category of groups with operators, from \cite{po} the category
$\cm(\ccl^1\gcl)$ of crossed modules in $\ccl^1\gcl$ is equivalent to the category
$Cat(\ccl^1\gcl)$ of internal categories in $\ccl^1\gcl$. On the other hand there is an
equivalence of categories between $Cat(\ccl^1\gcl)$ and $\ccl^2\gcl$, as explained for
instance in the proof of \cite[I--6, Proposition 1.2.3]{ell1}. It follows that
$\cm(\ccl^1\gcl)$ is equivalent to $\ccl^2\gcl$. The correspondences giving this
equivalence of categories can be easily made explicit from \cite{po} and \cite{ell1}. The
category $\ccl^2\gcl$ of cat$^2$-groups is also equivalent to the category $Crs^2$ of
crossed squares, \cite{lodb}. Hence there is an equivalence of categories between
$\cm(\ccl^1\gcl)$ and $Crs^2$. The correspondences giving this equivalence can be
described explicitly as follows:
\begin{lemma}\label{catone.lem1}\emph{ }

    a) Let $((H,d_{0},d_1),(H',d'_0,d'_1),\alpha)$ be an object of $\cm(\ccl^{1}\gcl)$. Let
    $(T,G,\mu)$ and $(T',G',\mu')$ be the crossed modules corresponding to
    $(H,d_{0},d_1)$ and $(H',d'_0,d'_1)$ respectively. Then the following is a crossed
    square
\begin{diagram}[s=2em]
    T & \rTo^{\mu} & G &&&& h:G\times T'\rw T \\
    \dTo^{\alpha_{|_{T\times 1}}} &&\dTo_{\alpha_{|_{1\times G}}} &&&&
    h(g,t')=(1,g)\;^{(t',1)}(1,g^{-1}) \\
    T' & \rTo_{\mu'}& G'
\end{diagram}
where $^{(t',1)}(1,g^{-1})$ is the crossed module action of $H'\cong T'\rt G'$ on $H\cong
T\rt G$.

b) Conversely, if
\begin{diagram}[s=2em]
    T & \rTo^\alpha & T' &&&& h:T'\times G\rw T\\
    \dTo^\mu && \dTo_{\mu'}\\
    G & \rTo_{\beta} & G'
\end{diagram}
is a crossed square, the corresponding object of $\cm(\ccl^1\gcl)$ is
\begin{equation*}
  (T\rt G,d_0,d_1)\xrw{(\alpha,\beta)}(T'\rt G',d'_0,d'_1)
\end{equation*}
where $d_0(t,g)=(1,g)$, $\;d_1(t,g)=(1,\mu(t)g)$, $\;d'_0(t',g')=(1,g')$,
$\;d'_1(t',g')=(1,\mu(t')g')$ for all $(t,g)\in T\rt G$, $\;(t',g')\in T'\rt G'$ and the
action of $T'\rt G'$ on $T\rt G$ is given by
\begin{equation*}
  ^{(t',g')}(t,g)=(\,^{t'}(\,^{g'}t)\,h\,(t',\,^{g'}g),\:^{g'}g)
\end{equation*}
for all $(t',g')\in T'\rt G'$, $\;(t,g)\in T\rt G$.
\end{lemma}
\prf It follows from the correspondences giving the equivalence of categories between
$\cm(\ccl^1\gcl)$ and $\ccl^2\gcl$ (see \cite{po} and \cite{ell1}) and between
$\ccl^2\gcl$ and $Crs^2$ (see \cite{lodb}).\enpr

%%%%%% New section %%%%%%%%%%%%%%%%%%%%%%%%%%%%%%%%%%%%%%%%%%%%%%%%%%%%%%%%%%%%%%%%%%%

\section[{Definition of the (co)homology and elementary properties}]
{Definition of the (co)homology and elementary properties} \label{ddfcm} In order to
define our (co)homology theory, we first introduce $\der$ and $\diff$ functors on the
category of crossed modules. We can do so in two equivalent ways, working directly in the
category $\cm$ or working in the equivalent category of cat$^1$-groups. We shall
illustrate both ways in some detail. While the first one may be slightly more
transparent, the reason we write explicitly the derivation functor in cat$^1$-groups is
that viewing our (co)homology theory as (co)homology of cat$^1$-groups will allow us to
apply in the next sections the interpretation in terms of extensions of the first and
second cotriple cohomology in a category of interest \cite{vale}.

Let $\Phi=\tgmu$ be a crossed module, $A$ an abelian group. The abelian crossed module
$(A,1,0)$ is a $\Phi$-module if and only if there is a split extension of crossed modules
\begin{equation}\label{ddfcm.eq1}
  (A,1,0)\remb(T',G',\mu')\;\pile{\suparle{}\\ \ronto}\;\tgmu
\end{equation}
Since the base group of the crossed module on the left is 1 the morphism
$(T'\!,G'\!,\mu'\!)\rw\tgmu$ is an isomorphism at the level of base groups and so we can
assume $G'=G$. Also observe that the section in (\ref{ddfcm.eq1}) gives, by conjugation
on $T'$, an action of $T$ on $A$, hence $T'\cong A\rt T$. So we can assume that the
sequence (\ref{ddfcm.eq1}) has the form
\begin{equation*}
  (A,1,0)\remb(A\rt T,G,\mu')\;\pile{\suparle{}\\ \ronto}\;\tgmu
\end{equation*}
where the morphisms on the left and on the right are the canonical inclusion and
projection respectively. The action of $G$ on $A\rt T$ induces an action of $G$ on $A$;
in fact, since $(\pr_T,\id_G)$ is a map of crossed modules, $\pr_T(\:^{g}(a,1))=1$ for
all $g\in G,\;\;a\in A$. Also, since the splitting $(i_T,\id_G)$ is a map of crossed
modules, $(0,\,^{g}t)=\:^{g}(0,t)$ for all $g\in G, \;\;t\in T$. It follows that, for all
$g\in G,\;\; t\in T$, $\;\widetilde{\mu}(a,t)=\mu(t)$,
$\;^{g}(a,t)=\:^{g}(a,1)\:^{g}(0,t)=(\:^{g}a,1)(0,\:^{g}t)$ and
$\widetilde{\mu}(a,t)=\mu(t)$. Requiring that the Peiffer identity holds for the crossed
module $(A\rtimes T,G,\widetilde{\mu})$ an easy calculation shows that for each $a,a'\in
A,\;\;t,t'\in T$
\begin{equation*}
  ^{\mu(t)}a'=a+\,^{t}a'-\,^{tt't^{-1}}a.
\end{equation*}
It follows that $a=\,^{t}a=\,^{\mu(t)}a$ for each $a\in A,\;\;t\in T$. Hence $A$ is a
$\p$-module, where $\p=G/\mu(T)$. We denote $A\rt \Phi=(A\times T,G,\widetilde{\mu})$.
From a general fact in algebraic categories,
\begin{equation*}
  \der(\Phi,(A,1,0))\cong\hund{\cm/\Phi}(\Phi,A\rt \Phi).
\end{equation*}
Since $A\rt \Phi=(A\times T,G,\widetilde{\mu})$ it is straightforward that
$\hund{\cm/\Phi}(\Phi,A\rt\Phi)\cong\hund{G}(T,A)$, where $\hund{G}(T,A)$ is the group of
$G$-equivariant homomorphisms from $T$ to $A$. In conclusion we obtain
\begin{equation}\label{ddfcm.eq2}
  \der(\tgmu,(A,1,0))\cong\hund{G}(T,A).
\end{equation}
The same procedure can be repeated in the equivalent category of cat$^1$-groups; let
$(T\rt G,d,s)$, $\;(A,0,0)$ be the cat$^1$-groups corresponding to the crossed modules
$\tgmu$ and (A,1,0) respectively. It is easy to see that $(A,0,0)$ is a $(T\rt
G,d,s)$-module if and only if there is a split extension of cat$^1$-groups
\begin{equation*}
  (A,0,0)\remb(A\rt(T\rt G),d',s')\;\pile{\suparle{}\\\ronto}\;(T\rt G,d,s)
\end{equation*}
where $d'(a,(t,g))=(0,d(t,g))$, $\;s'(a,(t,g))=(0,s(t,g))$, $\;a\in A$, $\;(t,g)\in T\rt
G$. Requiring that the identity $[\ker d',\ker s']=1$ holds, an easy calculation shows
that, for each $a,a'\in A$, $t,t'\in T$
\begin{equation*}
  a+\,^{(t,1)}a'-\,^{(t',\mu(t')^{-1})}a-a'=0.
\end{equation*}
It follows that $a=\,^{(t,1)}a=\,^{(1,\mu(t))}a$ for each $a\in A$, $t\in T$. Thus
$T\rtimes G$ acts on $A$, $T\times 1$ and $1\times \mu(T)$ act trivially on $A$, so $A$
is a $\p$-module, $\p=G/\mu(T)$.

From \cite{orz} a derivation from $(T\rtimes G,d,s)$ into $(A,0,0)$ is a map $D:T\rtimes
G\rw A$ which is a group derivation and such that it commutes with the 1-ary operations
$\omega\in\Omega_1\verb"\"\{-\}$. Hence $D(d(t,g))=D(s(t,g))=0$ for $(t,g)\in T\rtimes
G$, so that $D(1,g)=0$ for every $g\in G$.

In conclusion
\begin{equation}\label{ddfcm.eq3}
  \der((T\rt G,d,s),(A,0,0))\cong\{D\in\der(T\rt G,A)\;|\;D(1,G)=0\}.
\end{equation}
The two approaches are clearly equivalent. In fact there is an isomorphism
\begin{equation*}
  \alpha:\;\{D\in\der(T\rt G,A)\;|\;D(1,G)=0\}\rw\hund{G}(T,A)
\end{equation*}
given by $\alpha(D)(t)=D(t,1)$, as easily checked. This motivates the following
definition.
\begin{definition}\label{ddfcm.def1}
    Let $\Phi=\tgmu$ be a crossed module, $A$ an abelian group, $\p=G/\mu(T)$. We say
    that $\Phi$ acts on $A$ if $A$ is a $\p$-module. In this case we define
\begin{equation*}
  \der(\Phi,A)=\{D\in\der(T\rtimes G,A)\;|\;D(1,G)=0\}\cong\hund{G}(T,A).
\end{equation*}
where $\der(T\rtimes G,A)$ denotes group derivations from $T\rtimes G$ into $A$ and the
action of $T\rtimes G$ on $A$ is given by $(t,g)a=\:^{g\mu(T)}a$.
\end{definition}

Similarly if $\Phi$ acts on $A$ we have a contravariant functor
$\der(\mi,A):\cm/\Phi\rw\ab$ on the slice category; in fact, given an object
$\Phi'\rw\Phi$ of $\cm/\Phi$, the action of $\Phi$ on $A$ induces an action of $\Phi'$ on
$A$. When the context is clear, given an action of $\Phi$ on $A$ we will write
$\der(\Phi'\rw\Phi,A)$ as $\der(\Phi',A)$.

Given a crossed module $\Phi=\tgmu$ let $J_{T,G}$ be the ideal of $\zbb(T\rtimes G)$
generated by $\{(1,g)-(1,1)\;|\;1\neq g\in G\}$. Then
\begin{equation*}
  \der(\Phi,A)\cong\hund{\zbb(T\rtimes G)}\Big(\frac{\ifr_{\zbb(T\rtimes
  G)}}{J_{T,G}},A\Big)\cong\hund{\zbb \p}\Big(\zbb\p\otund{\zbb(T\rtimes G)}\frac{\ifr_{\zbb(T\rtimes
  G)}}{J_{T,G}},A\Big).
\end{equation*}
This motivates our next definition.

\begin{definition}\label{ddfcm.def2}
    Let $\Phi=\tgmu$ be a crossed module. Define
\begin{equation*}
  \diff\Phi=\zbb\p\otund{\zbb(T\rtimes G)}\frac{\ifr_{\zbb(T\rtimes G)}}{J_{T,G}}.
\end{equation*}
\end{definition}

For a crossed module $\Phi$ acting on the abelian group $A$, we denote
$\diff(\Phi,A)=A\otund{\zbb\p}\diff\Phi$. Similarly we have a covariant functor
$\diff(\mi,A):\cm/\Phi\rw\ab$.
\smallskip

The slice category $\cm/\Phi$ is tripleable over $\set/\mathcal{U}(\Phi)$ and we shall
denote by $\gbb$ the corresponding cotriple. We now consider the cotriple (co)homology of
$\Phi$ with coefficients in the $\Phi$-module $(A,1,0)$. This is the left (resp. right)
derived functor of the functor $\diff(\mi,A)$ (resp. $\der(\mi,A))$ on the slice category
$\cm/\Phi$ with respect to the cotriple $\gbb$.

\begin{definition}\label{defco.def1}
Let $\Phi=\tgmu$ be a crossed module, $A$ a \md{\p}. Define for each $n\geq 0$
\begin{equation*}
\begin{split}
   & D^n(\Phi,A)=H^n\der(\gbb_{\bullet}\Phi,A) \\
   & D_n(\Phi,A)=H_n\diff(\gbb_{\bullet}\Phi,A).
\end{split}
\end{equation*}
\end{definition}

The following are elementary properties of the (co)homology which can be deduced from
well known general facts about cotriple (co)homology in an algebraic category
\cite{bab1}. In what follows a projective crossed module means a projective object in the
category of crossed modules.

\begin{proposition}\label{defco.pro1}
Let $\Phi=\tgmu$ be a crossed module acting on the abelian group $A$. Then
\begin{align*}
a)& \; \;D^0(\Phi,A)\cong\der(\Phi,A),\qquad D_0(\Phi,A)\cong\diff(\Phi,A).\\
b)& \; \;\text{If } \Phi \text{ is a projective crossed module,}\\
 &\;\;D^n(\Phi,A)=0,\;\;D_n(\Phi,A)=0 \text{ for each } n>0. \\
c)& \; \;\text{Any short exact sequence of \md{\p}s } \;  \;0\rw A\rw A'\rw A''\rw 0 \;\;\text{induces long}\\
& \;\; \text{exact (co)homology sequences}\\
&  \cdots\rw D^n(\Phi,A)\rw D^n(\Phi,A')\rw D^n(\Phi,A'')\rw D^{n+1}(\Phi,A)\rw\cdots\\
 &   \cdots\rw D_n(\Phi,A)\rw D_n(\Phi,A')\rw D_n(\Phi,A'')\rw D_{n-1}(\Phi,A)\rw\cdots
\end{align*}
\end{proposition}

%%%%%%% New section %%%%%%%%%%%%%%%%%%%%%%%%%%%%%%%%%%%%%%%%%%%%%%%%%%%%%%%%%%%%%%
\section{(Co)homology of aspherical crossed modules}\label{ascro}
A crossed module $\Phi=\tgmu$ is called \emph{aspherical} when the map $\mu$ is
injective. The category of aspherical crossed modules is isomorphic to the category of
surjective group homomorphisms. Given a surjective group homomorphism $f:G\rw G'$ the
corresponding aspherical crossed module is $(\ker f,G,i)$ and will be denoted by
$\Phi_f$. If $\Phi_f$ acts on the abelian group $A$, then $A$ is a \md{\zbb G'} as well
as a \md{\zbb G} via $f$. We say in this case that $A$ is an \emph{\md{f}} and we denote
$\der(\Phi_f,A)=\der(f,A)$, $\;\diff \Phi_f=\diff f$.

\begin{lemma}\label{ascro.lem1}
    Let $f:G\rw G'$ be a surjective group homomorphism, $N=\ker f$, $\;A$ an \md{f}. Then
\begin{align*}
   &a)\quad\der(f,A)\cong \hund{\zbb G'}(N_{ab},A),\qquad \diff f\cong N_{ab}.\\
  &b)\quad\text{Suppose that there is a group homomorphism $f':G'\rw
  G$ with $ff'=\id$. Then}\\
  &\quad\;\;\;\text{there are short exact sequences}\\
  &\quad \qquad\qquad 0\rw\der(G',A)\rw\der(G,A)\rw\der(f,A)\rw 0\\
  &\quad \qquad\qquad 0\rw A\otund{\zbb G'}\diff f\rw
  A\otund{\zbb G}\ifr_G\rw A\otund{\zbb
  G'}\ifr_{G'}\rw 0.
\end{align*}
\end{lemma}

\prf\\ a) Let $\alpha:\der(f,A)\rw\hund{\zbb G'}(N_{ab},A)$ be defined by
\begin{equation*}
  \alpha(D)(n[N,N])=D(n,1),\qquad n\in N. \qquad\qquad
\end{equation*}
It is straightforward that $\alpha(D)$ is well defined; it is also a $\zbb
G'$-homomorphism since, for each $g'=f(g)\in G',\;n\in N$
\begin{equation*}
\begin{split}
   & \alpha(D)(g'\cdot n[N,N])=\alpha(D)(gng^{-1}[N,N])=D((1,g)(n,1)(1,g^{-1}))= \\
   & =(1,g)D(n,1)=g'\alpha(D)(n[N,N]).
\end{split}
\end{equation*}
Let $\beta:\hund{\zbb G'}(N_{ab},A)\rw\der(f,A)$ be defined by
\begin{equation*}
  \beta(\varphi)(n,g)=\varphi(n[N,N]),\qquad (n,g)\in N\rtimes G.
\end{equation*}
Then $\,\beta(\varphi)\,$ is a derivation and $\,\beta(\varphi)\,(1,G)\,=\,0\,$. We have
$\,\alpha\,\beta(\varphi)\,(n[N,N])=\beta(\varphi)(n,1)=\varphi(n[N,N])$ for each $n\in
N$; for each $D\in\der(f,A)$, $\;(n,g)\in N\rtimes G\,\,$ it is
$\;\beta\alpha(D)(n,g)=\alpha(D)(n[N,N])=D(n,1)=D(n,g)$. Thus $\alpha$ and $\beta$ are
inverse bijections and $\der(f,A)\cong\hund{\zbb G'}(N_{ab},A)$. Since this isomorphism
holds for each $\zbb G'$-module $A$, Yoneda Lemma implies $\diff f\cong N_{ab}$.
\smallskip

\noindent b) Let $\alpha:\der(G,A)\rw\der(f,A)$ be defined by $\alpha(D)(n,g)=D(n)$,
$\;D\in\der(G,A)$, $\;(n,g)\in N\rtimes G$. If $\;\xi\in\der(f,A)\;$, let
$\;D(g)=\xi(gf'f(g^{-1}),1), \;g\in G$. For each $g_1,g_2\in G$
\begin{equation*}
\begin{split}
   &  \;\xi(g_1g_2f'f(g_2^{-1})\,f'
    \,f(g_1^{-1}),1)\,=\,\xi((1,g_1)\,(g_2f'f(g_2^{-1}),1)\,(1,g_1^{-1})
    (g_1f'f(g_1^{-1}),1)=\\
   & =\xi(g_1f'f(g_1^{-1}),1)\,+\,(1,g_1)\;\xi(g_2f'f(g_2^{-1}),1).
\end{split}
\end{equation*}
This shows that $D$ is a derivation. Also, $\,(\alpha D)(n,g)=D(n)=\xi(n,1)=\xi(n,g)$ so
$\alpha$ is surjective. Exactness at the other terms is straightforward and b) follows
for the cohomology case.

From the well known exact sequence $N_{ab}\remb\zbb G'\otund{\zbb
G}\ifr_G\ronto\ifr_{G'}$ and from part a) we obtain the exact sequence
\begin{equation*}
    A\otund{\zbb G'}\diff f\supar{\beta'}A\otund{\zbb G}\ifr_G\supar{\alpha'}A\otund{\zbb
G'}\ifr_{G'}\rw 0
\end{equation*}
where $\;\beta'(a\otimes n[N,N])=a\otimes (n-1)$, $\;a\in A,\;n\in N$,
$\;\alpha'(a\otimes\sum_i b_i(g_i-1))=a\otimes\sum_i b_i(f(g_i)-1)$, $\;a\in A, \; b_i\in
\zbb$, $\;1\neq g_i\in G$. Let $\gamma:A\otund{\zbb G}\ifr_G\rw A\otund{\zbb G'}\diff f$
be defined by
\begin{equation*}
  \gamma(a\otimes\sideset{}{_i}\sum x_i(g_i-1))=\sideset{}{_i}\sum a\otimes x_i(f'f(g_i^{-1})g_i)[N,N]\qquad\qquad
\end{equation*}
for $a\in A$, $\;x_i\in\zbb$, $\; 1\neq g_i\in G$. Then $\gamma\beta'(a\otimes
n[N,N])=\gamma(a\otimes(n-1))=a\otimes n[N,N]$ for each $a\in A$, $n\in N$, so
$\gamma\beta'=\id$. It follows that $\ker\beta'=0$ and b) is proved. \enpr

\smallskip

We recall the notion of \emph{relative group (co)homology} in the sense of \cite{loda}.
Let $f:G\rw G'$ be a surjective group homomorphism and $A$ an \md{f}. Let $C^*(G,A)$ and
$C_*(G,A)$ be the standard (co)chain complexes for computing group (co)homology. For each
$n\geq 0 $ define
\begin{equation*}
\begin{split}
   & H^n(G',G;A)=H^n\mathrm{coker\,}(C^*(G',A)\remb C^*(G,A)) \\
   & H_n(G',G;A)=H_n{\ker}(C_*(G,A)\ronto C_*(G',A)).
\end{split}
\end{equation*}
Notice that this definition differs from the one in \cite{loda} by a dimension shift of
1.
\smallskip
\begin{theorem}\label{ascro.the1}
    Let $f:G\rw G'$ be a surjective group homomorphism, $A$ an $f$-module. Then
\begin{equation*}
  D^n(\Phi_f,A)\cong
  \begin{cases}
    \der(f,A) & n=0 \\
    H^{n+1}(G',G;A) & n>0
  \end{cases}
\quad
   D_n(\Phi_f,A)\cong
  \begin{cases}
  \diff(f,A) & n=0 \\
    H_{n+1}(G',G;A) & n>0.
  \end{cases}
\end{equation*}
\end{theorem}
\smallskip
\prf Let $\gbb_{\bullet}\Phi_f=(T_{\bullet},G_{\bullet},i_{\bullet})$,
$\;S_{\bullet}=G_{\bullet}/i_{\bullet}(T_{\bullet})$ and
$\;\psi_{\bullet}:G_{\bullet}\ronto S_{\bullet}$ be the quotient maps. Since $\Phi_f$ is
aspherical, from \cite{ccg} $T_{\bullet}\rw\ker f$, $\;G_{\bullet}\rw G$,
$\;S_{\bullet}\rw G'$ are free simplicial resolutions and there is a short exact sequence
of free simplicial groups $T_{\bullet}\remb G_{\bullet}\ronto S_{\bullet}$. Let $\bot$ be
the free cotriple on $\mathrm{Groups}$ and
$\bot_{\bullet}f:\bot_{\bullet}G\ronto\bot_{\bullet}G'$ the induced homomorphisms. Since,
for each $n$, $\,\psi_n$ and $\bot_nf$ have a section, from Lemma \ref{ascro.lem1} there
is a commutative diagram of cochain complexes
\begin{diagram}[s=1em,h=1.5em]
    0& \rTo & \der(S_{\bullet},A) & \rTo & \der(G_{\bullet},A) & \rTo &
    \der(\psi_{\bullet},A) & \rTo & 0\\
    && \dTo_{\wr} && \dTo_{\wr} && \dTo\\
    0& \rTo & \der(\bot_{\bullet}G',A) & \rTo & \der(\bot_{\bullet}G,A) & \rTo &
    \der(\bot_{\bullet}f,A) & \rTo & 0\,,
\end{diagram}
where $\sim$ are cochain homotopy equivalences. Taking the corresponding long exact
cohomology sequences in each row and applying the five Lemma we deduce that
\begin{equation*}
  H^n\der(\psi_{\bullet},A)\cong H^n\der(\bot_{\bullet}f,A)
\end{equation*}
for each $n\geq 0$. On the other hand there is a commutative diagram of cochain complexes
\begin{diagram}[s=1em,h=1.5em]
    0& \rTo & \der(\bot_{\bullet}G',A) & \rTo & \der(\bot_{\bullet}G,A) & \rTo &
    \der(\bot_{\bullet}f,A) & \rTo & 0\\
    && \dTo && \dTo && \dTo\\
    0& \rTo & C^*(G',A) & \rTo^{\alpha_{\bullet}} & C^*(G,A) & \rTo &
    \mathrm{coker\,}\alpha_{\bullet} & \rTo 0\,,
\end{diagram}
where $\der(\bot_{\bullet}G,A)\rw C^*(G,A)$ and $\der(\bot_{\bullet}G',A)\rw C^*(G',A)$
are the natural cochain maps of the Barr-Beck theory which induce isomorphisms in
cohomology (see \cite{bab2}). Taking the long exact cohomology sequence in each row of
the above diagram and applying the five Lemma we deduce that for each $n>0$
\begin{equation*}
  H^n\der(\bot_{\bullet}f,A)\cong
  H^{n+1}\mathrm{coker\,}\alpha_{\bullet}=H^{n+1}(G,G';A).
\end{equation*}
The argument for homology is similar.\enpr

%%%%%%% New section %%%%%%%%%%%%%%%%%%%%%%%%%%%%%%%%%%%%%%%%%%%%%%%%%%%%%%%%%%%%%
\section{Further properties of the (co)homology}\label{furpro}
\subsection{The relationship with CCG (co)homology.}
In the next proposition we establish the relationship between the (co)homology theory
defined in Section 3 and the one in \cite{ccg}. For a crossed module $\tgmu$ let
$\zeta\tgmu=T$.

\begin{proposition}\label{ascro.pro1}
    Let $\Phi=\tgmu$ be a crossed module, $A$ a trivial \md{\p}. Then for each $n\geq 0$
\begin{equation*}
  D^n(\Phi,A)\cong H^{n+1}_{CCG}(\tgmu,(A,1,0)).
\end{equation*}
Suppose further that $A\cong\zbb$. Then for each $n\geq 0$
\begin{equation*}
  D_n(\Phi,\zbb)\cong\zeta H^{CCG}_{n+1}\tgmu.
\end{equation*}
\end{proposition}

\prf Let $\gbb_{\bullet}\tgmu=(T_{\bullet},G_{\bullet},i_{\bullet})$,
$\;S_{\bullet}=G_{\bullet}/i_{\bullet}(T_{\bullet})$. Since the action of $\p$ on $A$ is
trivial, using Lemma \ref{ascro.lem1} and the well known isomorphism
$T_n/[G_n,T_n]\cong\zbb\otund{\zbb S_n}(T_n)_{ab}$, we obtain for each $n\geq 0$
\begin{equation*}
\begin{split}
     & H^{n+1}_{CCG}(\tgmu,(A,1,0))= H^n\hund{\ab\cm}\bigg(\Big(\frac{T_{\bullet}}{[G_{\bullet},T_{\bullet}]},
     \frac{G_{\bullet}}{[G_{\bullet},G_{\bullet}]},i_{\bullet}\Big),(A,1,0)\bigg)\cong\\
     & \cong H^n\hund{\zbb} \bigg(\frac{T_{\bullet}}{[G_{\bullet},T_{\bullet}]},A\bigg)\cong H^n\hund{\zbb}(\zbb\otund{\zbb \p}\zbb\p\otund{\zbb S_{\bullet}}(T_{\bullet})_{ab},A)\cong \\
     & \cong H^n\hund{\zbb \p}(\diff\gbb_{\bullet}\Phi,A)=D^n(\Phi,A).\\
     &  \\
     &  D_n(\Phi,\zbb)=H_n(\zbb\otund{\zbb \p}\diff\gbb_{\bullet}\Phi)\cong H_n(\zbb\otund{\zbb\p}
     \zbb\p\otund{\zbb S_{\bullet}}(T_{\bullet})_{ab})\cong H_n\bigg( \frac{T_{\bullet}}{[G_{\bullet},T_{\bullet}]}\bigg)\\
     & \cong H_n\bigg(\zeta\Big(
     \frac{T_{\bullet}}{[G_{\bullet},T_{\bullet}]},\frac{G_{\bullet}}{[G_{\bullet},G_{\bullet}]},
     i_{\bullet}\Big)\bigg)\cong\zeta H_n(T_{\bullet},G_{\bullet},i_{\bullet})_{ab}=\zeta H^{CCG}_{n+1}\tgmu.
\end{split}
\end{equation*}
\enpr

\subsection{Interpretation of the first and second cohomology group.}
Let $\Phi=\tgmu$ be a crossed module acting on the abelian group $A$. We shall need a
notion of singular and two-fold special extensions of $\tgmu$ by $(A,1,0)$.
\begin{definition}\label{furpro.def1b}
    Let $\Phi=\tgmu$ be a crossed module acting on an abelian group $A$.
\begin{description}
  \item[i)] A singular extension of $\,\tgmu$ by $\,(A,1,0)$ is a
  short exact sequence of crossed modules
\begin{equation}\label{furpro.eq1c}
  (A,1,0)\remb (T',G,\mu')\overset{(f,\id_G)}{\ronto}\tgmu
\end{equation}
such that the corresponding short exact sequence of cat$^1$-groups
\begin{equation}\label{furpro.eq1b}
  (A\times 1,0,0)\remb (T'\rt G,d',s')\overset{(f,\id_G)}{\ronto}(T\rt G,d,s)
\end{equation}
is a singular extension of $(T\rt G,d,s)$ by the $(T\rt G,d,s)$-module $(A\times 1,0,0)$
in the sense of categories of interest \cite{vale}.
  \item[ii)] A 2-fold special extension of $\,\tgmu$ by $(A,1,0)$ is an exact sequence of
  crossed modules
\begin{equation}\label{furpro.eq1d}
  (A,1,0)\overset{i}{\remb} (T'',G'',\mu'')\supar{(\alpha,\beta)}(T',G',\mu')\overset{(f,r)}{\ronto}\tgmu
\end{equation}
such that the corresponding exact sequence of cat$^1$-groups
\begin{equation}\label{furpro.eq2b}
  (A\times 1,0,0)\overset{i}{\remb} (T''\rt G'',d'',s'')\supar{(\alpha,\beta)}(T'\rt G',d',s')
  \overset{(f,r)}{\ronto}(T\rt G,d,s)
\end{equation}
is a 2-fold special extension of $(T\rt G,d,s)$ by the $(T\rt G,d,s)$-module $(A\times
1,0,0)$ in the sense of categories of interest \cite{vale}.
\end{description}
\end{definition}
We now give a more explicit characterization of singular and 2-fold special extensions of
$\tgmu$ by $(A,1,0)$.
\begin{lemma}\label{furpro.lem1b}
    Let $\Phi=\tgmu$ be a crossed module acting on the abelian group $A$.
\begin{description}
  \item[i)] A singular extension of $\,\tgmu$ by $\,(A,1,0)$ consists of a short exact
  sequence of crossed modules (\ref{furpro.eq1c}) such that if $f':T\rw T'$ is a set map with $ff'=\id_T$, it is
\begin{equation}\label{furpro.eq4b}
  f'(t)af'(t^{-1})=a, \qquad ^{g}a=\,^{[g]}a
\end{equation}
for all $g\in G$, $\;a\in A$, where $[g]=g\mu(T)\in \p$, $^{[g]}a$ is the given
$\p$-module action on $A$ and $^{g}a$ is given by the crossed module action of $G$ on
$T'$.
  \item[ii)] A 2-fold special extension of $\,\tgmu$ by $\,(A,1,0)$ consists of an exact
  sequence of crossed modules (\ref{furpro.eq1d}) where
\begin{equation}\label{furpro.eq6b}
  \begin{diagram}[s=2em]
    T'' & \rTo^\alpha & T' &&&& h:T'\times G''\rw T'' \\
    \dTo^{\mu''} && \dTo_\mu\\
    G'' & \rTo_\beta & G'
\end{diagram}
\end{equation}
is a crossed square and if $f':T\rw T'$, $\;r':G\rw G'$ are set maps with $ff'=\id_T$,
$\;rr'=\id_G$, then for all $g\in G$, $t\in T$, $a\in A$
\begin{equation}\label{furpro.eq7b}
  ^{r'(g)}a=\,^{[g]}a, \qquad ^{f'(t)}a=a.
\end{equation}
Here $[g]=g\mu(T)\in \p$, $^{[g]}a$ is the given $\p$-module action on $A$ while
$^{r'(g)}a$ (resp. $^{f'(t)}a$) is the  action of $G'$ (resp. $T'$) on $T''$ in the
crossed square (\ref{furpro.eq6b}).
\end{description}
\end{lemma}
\prf \emph{ }

(i) By definition (\ref{furpro.eq1c}) is a singular extension of crossed modules if and
only if the corresponding extension of cat$^1$-groups (\ref{furpro.eq1b}) is a singular
extension in the sense of categories of interest. By definition this means that the
induced action of $(T\rt G,d,s)$ on $(A\times 1,0,0)$ coincides with the given action,
that is
\begin{equation*}
  (f'(t),g)(a,1)(f'(t),g)^{-1}=\:^{(t,g)}(a,1).
\end{equation*}
An easy calculation shows this is equivalent to
\begin{equation*}
  (f'(t)\:^{g}af'(t^{-1}),1)=(\,^{[g]}a,1).
\end{equation*}
Hence for all $t\in T$, $g\in G$, $A\in A$,
\begin{equation}\label{furpro.eq8b}
  f'(t)\:^{g}a f'(t^{-1})=\:^{[g]}a.
\end{equation}
It is immediate to check that (\ref{furpro.eq8b}) is equivalent to (\ref{furpro.eq4b}).

(ii) By definition (\ref{furpro.eq1d}) is a 2-fold special extension of crossed modules
if and only if the corresponding extension of cat$^1$-groups (\ref{furpro.eq2b}) is a
2-fold special extension in the sense of categories of interest. By definition (see
\cite{vale}) this means that
\begin{description}
  \item[a)] $(A\times 1,0,0)$ is a $(T\rt G,d,s)$-module,
  \item[b)] $((T''\rt G'', d'', s''),(T'\rt G',d',s'), (\alpha,\beta))$ is a crossed module in the category of
cat$^1$-groups,
  \item[c)] $(A\times 1,0,0)\supar{i}(T''\rt G'',d'',s'')$ is a morphism of $(T'\rt
G',d',s')$-structures, where $(T'\rt G',d',s')$ acts on $(A\times 1,0,0)$ via $(f,r)$.
\end{description}

Condition b) and Lemma \ref{catone.lem1} imply that (\ref{furpro.eq6b}) is a crossed
square, and the crossed module action of $T'\rt G'$ on $T''\rt G''$ is given by
\begin{equation}\label{furpro.eq9b}
  ^{(t',g')}(t'',g'')=(\,^{t'}(\,^{g'}t'')h(t',\,^{g'}g''),\,^{g'}g'').
\end{equation}
It easily checked that condition c) is equivalent to requiring that the induced action of
$(T\rt G,d,s)$ on $(A\times 1,0,0)$ given by $^{(f'(t),r'(g))}(a,1)$ coincides with the
given action which is $(\,^{[g]}a,1)$. Hence by (\ref{furpro.eq9b}) we obtain
\begin{equation*}
  (^{f'(t)}( ^{r'(g)}a)h(f'(t),1),1)=(\,^{[g]}a,1)
\end{equation*}
for all $t\in T$, $g\in G$ , $a\in A$. From the axioms of crossed squares \cite{lodb}
$h(f'(t),1)=1$, hence the above is equivalent to
\begin{equation}\label{furpro.eq10b}
  ^{f'(t)}(\,^{r'(g)}a)=\:^{[g]}a
\end{equation}
for all $t\in T$, $g\in G$ , $a\in A$. It is straightforward that (\ref{furpro.eq10b}) is
equivalent to (\ref{furpro.eq7b}).\enpr

\medskip

Two singular extensions of $\tgmu$ by $(A,1,0)$,
\begin{equation*}
  (A,1,0)\remb(T'_i,G,\mu'_i)\ronto\tgmu
\end{equation*}
$i=1,2$ are congruent if there is a morphism of crossed modules
$\psi:(T'_1,G,\mu'_1)\rw(T'_2,G,\mu'_2)$ such that the following diagram commutes
\begin{diagram}[s=2em]
    (A,1,0) & \rEmbed & (T'_1,G,\mu'_1) & \rOnto & \tgmu\\
    \dEq && \dTo_\psi && \dEq\\
    (A,1,0) & \rEmbed & (T'_2,G,\mu'_2)& \rOnto & \tgmu .
\end{diagram}
It follows that $\psi$ is an isomorphism. Hence congruence defines an equivalence
relation on the set of singular extensions of $\Phi$ by $(A,1,0)$ and we can consider the
set of equivalence classes $\mathcal{E}^1(\Phi,A)$. This is an abelian group with Baer
sum, the zero element being the class $[(A,1,0)\remb A\rt \Phi \ronto \Phi]$.

Two 2-fold special extensions of $\Phi=\tgmu$ by $(A,1,0)$ are related if there is a
morphism
\begin{diagram}[s=2em]
    (A,1,0) & \rEmbed & (T''_1,G''_1,\mu'_1) &\rTo & (T'_1,G'_1,\mu'_1) & \rOnto & \tgmu\\
    \dEq && \dTo^\alpha && \dTo_\beta && \dEq\\
    (A,1,0) & \rEmbed & (T''_2,G''_2,\mu'_2) & \rTo & (T'_2,G'_2,\mu'_2) & \rOnto & \tgmu
\end{diagram}
such that $(\alpha,\beta)$ is a morphism of crossed squares. This relation generates an
equivalence relation and we denote by $\mathcal{E}^2(\Phi,A)$ the set of equivalence
classes of 2-fold special extensions of $\Phi$ by $(A,1,0)$. This is in fact an abelian
group with Baer sum \cite{vale}.
\begin{proposition}\label{furpro.pro1}
    Let $\Phi$ be a crossed module acting on the abelian group $A$. There are
    isomorphisms of abelian groups
\begin{equation*}
\begin{split}
  D^1(\Phi,A)\cong\mathcal{E}^1(\Phi,A)\: \\
  D^2(\Phi,A)\cong\mathcal{E}^2(\Phi,A).
\end{split}
\end{equation*}
\end{proposition}
\prf From Section 2, the cohomology $D^*(\Phi,A)$ is cotriple cohomology in the category
of interest $\cm$ with coefficients in the \md{\Phi} $(A,1,0)$. The interpretation in
terms of extensions of the first and second cotriple cohomology in any category of
interest can be found for example in \cite{vale}. The result is thus a direct
specialization of \cite[Theorem 2.1.3, Proposition 2.1.5, Theorem 2.2.3, Proposition
2.2.4 ]{vale}. \enpr

\medskip

If $f:G\rw G'$ is a surjective group homomorphism and $A$ is an \md{f}, from Theorem
\ref{ascro.the1} and Proposition \ref{furpro.pro1} we deduce that
\begin{equation*}
  H^2(G',G;A)\cong\mathcal{E}^1(\Phi_f,A), \quad H^3(G',G;A)\cong\mathcal{E}^2(\Phi_f,A).
\end{equation*}
We observe that the first of these isomorphisms recovers a result of \cite{loda}. We
notice in fact that the group $\mathcal{E}^1(\Phi_f,A)$ is isomorphic to the group of
relative extensions of $(G',G)$ by $A$; these consist of exact sequences of groups
\begin{equation*}
  0\rw A\rw M \supar{\mu} G \supar{f} G'\rw 1
\end{equation*}
such that $\mu$ is a crossed module and the induced action of $G'$ on $A$ coincides with
the given one. A congruence of relative extensions is a commutative diagram
\begin{diagram}[h=1.4em,s=1.5em]
    0 & \rw & A & \rw & M & \supar{\mu} & G & \rw & G' & \rw & 1\\
    && \dEq && \dTo^\psi && \dEq && \dEq\\
    0 & \rw & A & \rw & Q & \supar{\mu'} & G & \rw & G' & \rw & 1
\end{diagram}
such that $(\psi,\id_G)$ is a morphism of crossed modules. Let $\mathcal{E}xt(G',G;A)$ be
the set of equivalence classes, made into an abelian group as in \cite{loda}. There is a
map of abelian groups $\alpha:\mathcal{E}xt(G',G;A)\rw\mathcal{E}^1(\Phi,A)$
\begin{equation*}
  \alpha\;[0\rw A\rw M \supar{\mu} G \supar{f} G'\rw 1]=[(A,1,0)\remb(M,G,\mu)\ronto(\ker
  f,G,i)].
\end{equation*}
It is immediate to check that $\alpha$ is well defined and that it is a bijection.

The interpretation of $H^3(G',G;A)$ in terms of equivalence classes of 2-fold special
extensions of $(\ker f,G,i)$ by $(A,1,0)$ does not seem to have been given in the
literature as far as the author knows. Notice that the identification of relative group
cohomology with cotriple cohomology of a crossed module also allows to give a simplicial
interpretation of $H^n(G,G';A)$ for any $n$ by direct application of the results of
\cite{dusk}.

\subsection{Universal coefficient formulae.} We shall establish universal coefficient
formulae for the (co)homology of a crossed module $\Phi$ acting trivially on an abelian
group $A$.

\begin{theorem}\label{furpro.the1}
    Let $\Phi=\tgmu$ be a crossed module acting trivially on an abelian group $A$. Then
    there are short exact sequences
    \begin{equation*}
  0\rw\ext^1_\zbb (D_{n-1}(\Phi,\zbb),A)\rw
  D^n(\Phi,A)\rw\homz(D_n(\Phi,\zbb),A)\rw 0
\end{equation*}
\begin{equation*}
  0\rw D_n(\Phi,\zbb)\otund{\zbb}A\rw
  D_n(\Phi,A)\rw\tor^\zbb_1(D_{n-1}(\Phi,\zbb),A)\rw 0.
\end{equation*}
\end{theorem}

\prf Since the action of $\Phi$ on $A$ is trivial, from Proposition \ref{ascro.pro1} we
know that $D^n(\Phi,A)\cong H^{n+1}(\Phi,(A,1,0))$ and $D_n(\Phi,\zbb)\cong\zeta
H^{CCG}_{n+1}(\Phi)$.

Since $\zbb$ has global dimension 1, there exists a resolution $A\rw I^\bullet$ of $A$ by
injective $\zbb$-modules with $I^m=0$ for $m\geq 2$; then $(A,1,0)\rw(I^\bullet,1,0)$ is
an injective resolution of $(A,1,0)$ which satisfies the hypotheses of \cite[Theorem 18
(iv) ]{ccg}. The cohomology universal coefficient sequence then follows from
\cite[Theorem 18 (iv)]{ccg}.

For the homology case, denote $\Phi_{\bullet}=\gbb_{\bullet}\Phi$ and let
$\varphi_{\bullet\bullet}$ be the double complex of abelian groups
\begin{equation*}
  \varphi_{\bullet\bullet}=P_{\bullet}\otund{\zbb \p}\diff\Phi_{\bullet}
\end{equation*}
where $P_{\bullet}\rw A$ is a projective $\zbb$-resolution of $A$ with $P_n=0$ for $n>1$
(such resolution exists since $\zbb$ has global dimension 1). For any crossed module
$\Phi$ acting trivially on $A$, since
  $\zbb\otund{\zbb\p}\diff\Phi\cong D_0(\Phi,\zbb)\cong\zeta H^{CCG}_{1}\Phi\cong\zeta\Phi_{ab}$
we have
\begin{equation*}
  A\otund{\zbb\p}\diff\Phi\cong A\otund{\zbb}\zeta\Phi_{ab}.
\end{equation*}
In particular
\begin{equation*}
  P_{\bullet}\otund{\zbb \p}\diff\Phi_{\bullet}\cong
  P_{\bullet}\otmz\zeta(\Phi_{\bullet})_\mathrm{ab}.
\end{equation*}
Since $D^1(\Phi_{q},A)=0$ for each $q$, from the cohomology universal coefficient
sequence and from Proposition \ref{ascro.pro1} we obtain
$\mathrm{Ext}^1_\zbb(\zeta(\Phi_q)_{ab},A)=0$ for every abelian group $A$. It follows
that $\zeta(\Phi_q)_{ab}$ is a projective \md{\zbb}. Therefore
\begin{align*}
   & H^v_p(P_{\bullet}\otund{\zbb \p}\diff \Phi_q)\cong H^v_p(P_{\bullet}\otmz\zeta
   (\Phi_q)_\mathrm{ab})\cong \\
   & \cong \tor_p^\zbb(\zeta(\Phi_q)_\mathrm{ab},A)\cong
  \begin{cases}
    A\otund{\zbb\p}\diff \Phi_q & p=0, \\
    0 & p>0.
  \end{cases}
\end{align*}
Taking homology again we obtain
\begin{equation*}
  H^h_q\,H^v_p(P_{\bullet}\otund{\zbb \p}\diff\Phi_{\bullet})\cong
  \begin{cases}
    D_q(\Phi,A) & p=0, \\
    0 & p>0.
  \end{cases}
\end{equation*}
Therefore the spectral sequence
\begin{equation*}
  H^h_q\,H^v_p(\varphi_{\bullet\bullet})\Rightarrow H_{p+q}\,\mathrm{Tot}\,\varphi_{\bullet\bullet}
\end{equation*}
collapses, giving $H_n\,\mathrm{Tot}\,\varphi_{\bullet\bullet}=D_n(\Phi,A)\quad n\geq 0$.
Consider the second spectral sequence
\begin{equation*}
  H^v_q\,H^h_p(\varphi_{\bullet\bullet})\Rightarrow
  H_{p+q}\mathrm{Tot}\,\varphi_{\bullet\bullet}.
\end{equation*}
Fixing $p$ and taking homology we have:
\begin{align*}
    H^h_q(P_p\otund{\zbb \p}\diff \Phi_{\bullet})&\cong H^h_q(P_p\otmz\zeta(\Phi_{\bullet})_\mathrm{ab})\cong
   P_p\otmz H^h_q\zeta(\Phi_{\bullet})_\mathrm{ab}\cong\\
   &\cong P_p\otmz\zeta\,H^{CCG}_{q+1}(\Phi)\cong
  \begin{cases}
    P_p\otmz D_q(\Phi,\zbb) & p=0,1, \\
    0 & p>1.
  \end{cases}
\end{align*}
Taking homology again:
\begin{equation*}
  H^v_p\,H^h_q(P_{\bullet}\otund{\zbb \p}\diff\Phi_{\bullet})=
  \begin{cases}
    \tor^\zbb_p(D_q(\Phi,\zbb),A) & p=0,1, \\
    0 & p>1.
  \end{cases}
\end{equation*}
So we obtain a universal coefficient spectral sequence
\begin{equation*}
  E^2_{pq} \Rightarrow D_{p+q}(\Phi,A)
\end{equation*}
which has $E^2_{pq}=0$ for $p\neq 0,1$.

Therefore there are short exact sequences $0\rw E_{0n}^2\rw
H_n\mathrm{Tot}\,\varphi_{\bullet\bullet}\rw E^2_{1,n-1}\rw 0$, i.e.
\begin{equation*}
  0\rw D_n(\Phi,\zbb)\otmz A\rw D_n(\Phi,A)\rw \tor^\zbb_1(D_{n-1}(\Phi,\zbb),A)\rw 0.
\end{equation*}
\enpr

%%%%% New section %%%%%%%%%%%%%%%%%%%%%%%%%%%%%%%%%%%%%%%%%%%%%%%%%%%%%%%%%%%%%%%
\section{The relationship with the (co)homology of the classifying space.}\label{class}
In \cite{elli} the (co)homology of a crossed module $\Phi$ with coefficients in a
$\p$-module $A$ is defined as the (co)homology of the classifying space $B(\Phi)$ of the
crossed module with coefficients in the local system corresponding to $A$. We recall the
algebraic description of this (co)homology. This is a special case of a more general
construction, which is well known.

If $G_*$ is a simplicial group and $A$ is a $\pi_0(G_*)$-module, since $\p BG_*\cong\pi_0
G_*,\; $ $A$ is a local system on the classifying space $BG_*$ of $G_*$. There is an
algebraic description of the (co)homology $H_*(BG_*,A)$ and $H^*(BG_*,A)$. If $C_*(G,A)$
(resp. $C^*(G,A)$) is the standard chain (resp. cochain) complex for computing group
homology (resp. cohomology) then there are isomorphisms (see for instance \cite[Lemma
5.1]{pira}):
\begin{equation}\label{class.eq1a}
    H_*(\mathrm{Tot}(C_*(G_*,A))\cong H_*(BG_*,A),\qquad H^*(\mathrm{Tot}(C^*(G_*,A))\cong H^*(BG_*,A).
\end{equation}
If $\,G_*\rw H_*$ is a map of simplicial groups which is a weak equivalence, that is such
that it induces isomorphisms of homotopy groups, and $A$ a \md{\pi_0(G_*)}, then the
induced maps $H_*(BG_*,A)\rw H_*(BH_*,A)$ and $H^*(BG_*,A)\rw H^*(BH_*,A)$ are
isomorphisms.

Let $N_*^{-1}\tgmu$ be the simplicial group whose Moore complex has length 1
corresponding to the crossed module $\Phi=\tgmu$ as in Section 1. Taking
$G_*=N_*^{-1}\tgmu$ in (\ref{class.eq1a}) we obtain the algebraic description of the
(co)homology of the classifying space of the crossed module. If a morphism of crossed
modules $\alpha:\Phi\rw\Phi'$ is a weak equivalence then the (co)homology groups of the
classifying spaces of $\Phi$ and $\Phi'$ with coefficients in a $\p$-module $A$ are
isomorphic.
\smallskip

Our main result in this section is that the (co)homology of crossed modules defined in
Section 2 is related by a long exact sequence to the (co)homology of the classifying
space of the crossed module. In proving this result, we also give a simplicial
description of the (co)homology of the classifying space. An application of this will be
given in Section 7.

In the second part of this section we give an alternative description of the (co)homology
$D_*(\Phi,A)$ and $D^*(\Phi,A)$ without using cotriples.

\begin{theorem}\label{class.the1b}
    Let $\Phi=\tgmu$ be a crossed module, $A$ a $\p$-module. Let
    $\gbb_{\bullet}\Phi=(T_{\bullet},G_{\bullet},\mu_{\bullet})$ and
    $S_{\bullet}\cong G_{\bullet}/\mu_{\bullet}(T_{\bullet})$. Then

\medskip
\begin{description}
  \item[i)] $BS_{\bullet}$ and $B\Phi$ are weakly homotopy equivalent.
  \item[ii)]
  \begin{equation*}
 H_n\diff(S_{\bullet},A)\cong
  \begin{cases}
    H_{n+1}(B\Phi,A) &\quad n>0, \\
    A\otund{\zbb \p}\ifr_{\p} &\quad n=0.
  \end{cases}
  \end{equation*}
  \begin{equation*}
   \text{\emph{ }}\quad\: H^n\der(S_{\bullet},A)\cong
  \begin{cases}
    H^{n+1}(B\Phi,A) &\quad n>0, \\
    \der(\p,A) &\quad n=0,
  \end{cases}\quad\;\;
\end{equation*}
  \item[iii)] There are long exact (co)homology sequences
  \begin{equation*}
\begin{split}
  \cdots&\rw D_n(\Phi,A)\rw H_{n+1}(G,A)\rw H_{n+1}(B(T,G,\mu),A)\rw\\
  &\rw D_{n-1}(\Phi,A)\rw\cdots\rw A\otund{\zbb \p}\diff \Phi\rw A \otund{\zbb \p}\diff
  G\rw \\
  & \rw A \otund{\zbb \p}\diff \p\rw 0
  \end{split}
\end{equation*}
\begin{equation*}
\begin{split}
  0 & \rw\der(\p,A)\rw\der(G,A)\rw\der(\Phi,A)\rw H^2(B(T,G,\mu),A)\rw\\
    & \rw H^2(G,A)\rw D^1(\Phi,A)\rw\cdots
  \end{split}
\end{equation*}
\end{description}

\end{theorem}
\prf

 (i) Let $s\cm$ be the category of simplicial crossed modules, SimplSet (resp.
$\biset$) the category of simplicial (resp. bisimplicial) sets. Let $N_*^{-1}:\cm\rw
{\mathbf{SG}_{\leq 1 }}$ be as in \S \ref{crossm} and let $\mathcal{N}:\mathbf{SG}\rw
\text{SimplSet}$ be the functor associating to a simplicial group the diagonal of the
bisimplicial set obtained by forming the nerve of the group in each dimension of the
simplicial group. The composite $\mathcal{N}\circ N^{-1}_*$ is a functor $\cm\rw
\text{SimplSet}$. By definition, the classifying space of a crossed module $\Phi$ is the
geometric realization of the simplicial set $\mathcal{N} N^{-1}_*(\Phi)$. Given a
simplicial crossed module, we can apply $\mathcal{N}\circ N^{-1}_*$ in each dimension to
obtain a bisimplicial set. Hence we have a functor $F:s\cm\rw\biset$.

Consider in particular the simplicial crossed modules
$\gbb_{\bullet}\Phi=(T_{\bullet},G_{\bullet},\mu_{\bullet})$ and $\Phi=\tgmu$ (the second
is a constant simplicial crossed module).

We claim that $F(\gbb_{\bullet}\Phi)\rw F(\Phi)$ is a pointwise weak equivalence of
bisimplicial sets in the sense that all the maps $F(\gbb_{\bullet}\Phi)_{*m}\rw
F(\Phi)_{*m} $ are weak equivalences of simplicial sets. In fact, denoting by $\{\cdot\}$
the one point set
\begin{equation*}
\begin{split}
    (F\gbb_{\bullet}\Phi)_{nm}=(\mathcal{N}N_*^{-1}(T_n,G_n,\mu_n))_m=&
  \begin{cases}
    \{\cdot\} & m=0, \\
    T_n^{m^2}\times G_n^m & m>0.
  \end{cases}
    \\
    (F\Phi)_{nm}=(\mathcal{N}N^{-1}_*\tgmu)_m=&
  \begin{cases}
   \{\cdot\} & m=0, \\
   T^{m^2}\times G^m & m>0.
  \end{cases}
\end{split}
\end{equation*}
By the properties of the cotriple resolution  $\gbb_{\bullet}\Phi$ \cite{ccg},
$G_{\bullet}\rw G$ and $T_{\bullet}\rw T$ are free simplicial resolutions of groups,
therefore $T_{\bullet}^{m^2}\times G_{\bullet}^m\rw T^{m^2}\times G^m$ is a weak
equivalence in SimplSet, hence for each $m$, $\;F(\gbb_{\bullet}\Phi)_{*m}\rw
F(\Phi)_{*m}$ is a weak equivalence of simplicial sets, that is $F(\gbb_{\bullet}\Phi)\rw
F(\Phi)$ is a pointwise weak equivalence in $\biset$. It is proved in \cite[Ch. IV,
Proposition 1.7]{goer} that if $f:X\rw Y$ is a pointwise weak equivalence of bisimplicial
sets, in the sense that all the maps $f:X_m\rw Y_m $ are weak equivalences of simplicial
sets, then the induced map $f_*:\mathrm{diag\,}(X)\rw\mathrm{diag\,}(Y)$ of associated
diagonal simplicial sets is a weak equivalence. It follows that
$\mathrm{diag\,}F(\gbb_{\bullet}\Phi)\rw \mathrm{diag\,}F(\Phi)$ is a weak equivalence in
SimplSet, so that the respective geometric realizations
$|\mathrm{diag\,}F(\gbb_{\bullet}\Phi)|$ and $|\mathrm{diag\,}F(\Phi)|$ are weakly
homotopy equivalent.

Recall (see for instance \cite[p. 94]{quil}) that the geometric realization
$|\mathrm{diag\,}X|$ of the diagonal of a bisimplicial set $X$ is homeomorphic to the
geometric realization of the simplicial space obtained by taking the geometric
realization in vertical directions and is also homeomorphic to the geometric realization
of the simplicial space obtained by taking the geometric realization in the horizontal
directions. Hence $|\mathrm{diag\,}F(\gbb_{\bullet}\Phi)|$ is homeomorphic to the
geometric realization of the simplicial space
$\{|\mathcal{N}N^{-1}_*(T_n,\!G_n,\!\mu_n)|\!\}\!=\!\{B(T_n,\!G_n,\!\mu_n)\!\}$; but
since each crossed module $(T_n,\!G_n,\!\mu_n)$ is aspherical, $B(T_n,G_n,\mu_n)\cong
B(1,S_n,i)\cong BS_n$. The geometric realization of the simplicial space $\{BS_n\}\,$ is
homeomorphic to the geometric realization of the simplicial group
$S_{\bullet},BS_{\bullet}$. So in conclusion
$\;|\mathrm{diag\,}F(\gbb_{\bullet}\Phi)|\cong BS_{\bullet}$.

On the other hand clearly $|\mathrm{diag\,}F(\Phi)|=B\Phi$, so that $B\Phi$ and
$BS_{\bullet}$ are weakly homotopy equivalent, proving (i).

\smallskip

(ii) It follows from (i) that for each $n\geq 0$, $H_n(B\Phi,A)\cong
H_n(BS_{\bullet},A)$. On the other hand we observe that for each $n>1$
\begin{equation}\label{class.eq1b}
  H_n(BS_{\bullet},A)\cong H_{n-1}\diff(S_{\bullet},A)
\end{equation}
In fact, we can compute $H_n(BS_{\bullet},A)$ as indicated in (\ref{class.eq1a}). Since
$S_m$ is free, $0\rw \ifr_{S_m}\rw \zbb S_m \rw \zbb \rw 0$ is a free resolution of the
trivial \md{S_m} $\zbb$, hence $C_{\bullet}(S_m,A)$ is the complex $0\rw A\otund{\zbb
S_m}\ifr_{S_m}\supar{}A\rw 0$. It follows that $H_*(BS_{\bullet},A)$ is the total
homology of the bicomplex $\psi_{_{\bullet\bullet}}$
\begin{diagram}[h=1.5em,s=1.6em]
    && 0 && 0 && 0 && 0 &&\\
    && \dTo && \dTo && \dTo && \dTo && \\
    \cdots&\rTo & A\otund{\zbb S_3}\ifr_{S_3} & \rTo & A\otund{\zbb S_2}\ifr_{S_2} &
    \rTo & A\otund{\zbb S_1}\ifr_{S_1} & \rTo & A\otund{\zbb S_0}\ifr_{S_0} & \rTo & 0\\
    && \dTo && \dTo && \dTo && \dTo\\
    \cdots&\rTo & A & \rTo_{\id} & A & \rTo_{\id} & A & \rTo_{\id} & A & \rTo & 0 \\
    && \dTo && \dTo && \dTo && \dTo && \\
    && 0 && 0 && 0  && 0 &&.
\end{diagram}
The spectral sequence of this double complex is
\begin{equation*}
  E^1_{pq}=H^h_q(\psi_{*p})=
  \begin{cases}
    0 & p=0,\; q>0 \text{ or } p>1, \\
    A & p=0,\; q=0,\\
    H_q(A\otund{\zbb S_{\bullet}}\ifr_{S_{\bullet}}) & p=1,\; q\geq 0.
  \end{cases}
\end{equation*}
Therefore
\begin{equation*}
  E^2_{pq}=
  \begin{cases}
    0 & \text{for }  p\neq 0,1 \text{ or } p=0,\; q>0\\
    H_q(A\otund{\zbb S_{\bullet}}\ifr_{S_{\bullet}}) & p=1,\; q>0.
  \end{cases}
\end{equation*}
So there are short exact sequences $0\rw E^2_{0n}\rw
H_n\mathrm{Tot}\psi_{\bullet\bullet}\rw E^2_{1,n-1}\rw 0$ and since $E^2_{0n}=0$ for
$n>0$ and $E^2_{1,n-1}=H_{n-1}(A\otund{\zbb S_{\bullet}}\ifr_{S_{\bullet}})$ for $n>1$ we
deduce
\begin{equation*}
  H_n(BS_{\bullet},A)\cong H_n\mathrm{Tot}\psi_{\bullet\bullet}\cong
  E^2_{1,n-1}=H_{n-1}(A\otund{\zbb S_{\bullet}}\ifr_{S_{\bullet}})
\end{equation*}
for $n>1$, which is (\ref{class.eq1b}). It follows that $H_n\diff(S_{\bullet},A)\cong
H_{n+1}(B\Phi,A)$ for $n\geq 1$.

It remains to prove that $H_0\diff(S_{\bullet},A)\cong A\otund{\zbb \p}\ifr_{\p}$.
Consider the following diagram:
\smallskip
\begin{diagram}[h=1.2em,s=1.5em]
    0 & \rTo & A\otund{\zbb \p}\diff \gbb_1\Phi & \rTo & A\otund{\zbb \p}\diff G_1 &
    \rTo & A\otund{\zbb \p}\diff S_1 & \rTo & 0 \\
    && \dTo && \dTo && \dTo \\
    0 & \rTo & A\otund{\zbb \p}\diff \gbb_0\Phi & \rTo & A\otund{\zbb \p}\diff G_0 &
    \rTo & A\otund{\zbb \p}\diff S_0 & \rTo & 0 \\
    && \dTo && \dTo && \dTo\\
    &  & A\otund{\zbb \p}\diff \Phi & \rTo & A\otund{\zbb \p}\diff G &
    \rTo & A\otund{\zbb \p}\diff \p & \rTo & 0 \\
    && \dTo && \dTo && \dTo \\
    && 0 && 0 && 0 && .
\end{diagram}
The first two rows from the top are exact by Lemma \ref{ascro.lem1}. We claim that the
bottom row is also exact. In fact, it is straightforward to check that there is an exact
sequence
\begin{equation}\label{class.eq24}
  0\rw\der(\p,A)\supar{\alpha}\der(G,A)\supar{\beta}\der(\Phi,A)
\end{equation}
where $\alpha(D)(g)=D(g\mu(T)),\;\;g\in G,\;\;D\in\der(\p,A)$ and $\beta(D)=Dd_1-Dd_0$,
$\;d_1(t,g)=\mu(t)g$, $\;d_0(t,g)=g$, $\;(t,g)\in T\rtimes G$. Consider the map
$\gamma:\diff\Phi\rw\zbb\p\otund{\zbb G}\ifr_G$ defined by
\begin{equation*}
  \gamma\Big(x\otimes\Big[\sideset{}{_i}\sum a_i(y_i-e)\Big]\Big)=x\otimes\Big(\sideset{}{_i}\sum
  a_i(d_1(y_i)-1)-\sideset{}{_i}\sum a_i(d_0(y_i)-1)\Big)
\end{equation*}
$x\in\zbb\p$, $\;e\neq y_i\in T\rtimes G$, $\;a_i\in\zbb$. By left exactness of
$\hund{\zbb\p}(\mi,A)$ we obtain an exact sequence
\begin{equation*}
  0\rw\hund{\zbb\p}(\mathrm{coker\,}\gamma,A)\rw\der(G,A)\rw\der(\Phi,A).
\end{equation*}
Since $\hund{\zbb\p}(\gamma,A)=\beta$, (\ref{class.eq24}) implies that
$\hund{\zbb\p}(\mathrm{coker\,}\gamma,A)\cong \hund{\zbb\p}(\ifr_{\p},A)$ for every
\md{\zbb\p} $A$; hence $\mathrm{coker\,}\gamma\cong\ifr_{\p}$ and we have the exact
sequence
\begin{equation*}
  \diff\Phi\rw\zbb\p\otund{\zbb G}\ifr_G\rw\ifr_{\p}\rw 0.
\end{equation*}
From right exactness of $A\otund{\zbb\p}\mi$ the claim follows. Since the first two
columns from the left of the diagram are also exact, it follows from an easy diagram
chasing argument that the third column is also exact. Therefore $H_0(A\otund{\zbb
\p}\diff S_{\bullet})\cong A\otund{\zbb \p}\diff \p$. The argument for cohomology is
similar.

\smallskip

(iii) From \cite{ccg} $\;T_{\bullet}\rw T$ and $G_{\bullet}\rw G$ are free simplicial
resolutions and there is a short sequence of free simplicial groups $T_{\bullet}\remb
G_{\bullet}\overset{\psi_{\bullet}}{\ronto} S_{\bullet}$. From Lemma \ref{ascro.lem1} we
have short exact sequences of (co)chain complexes
\begin{equation*}
\begin{split}
   & 0\rw A\otund{\zbb \p}\diff\psi_{\bullet}\rw A\otund{\zbb \p}\diff G_{\bullet}\rw
   A\otund{\zbb \p}\diff S_{\bullet}\rw 0\\
   & 0\rw\der(S_{\bullet},A)\rw\der(G_{\bullet},A)\rw\der(\psi_{\bullet},A)\rw 0.
\end{split}
\end{equation*}
Taking the corresponding long exact (co)homology sequences and using (ii) we obtain
\begin{equation}\label{class.eq16}
\begin{split}
   & \cdots\rw D_n(\Phi,A)\rw H_{n+1}(G,A)\rw H_n\diff(S_{\bullet},A)\rw D_{n-1}(\Phi,A)\rw \cdots\\
   & \rw A\otund{\zbb \p}\diff\Phi\rw A\otund{\zbb \p}\diff G\rw
   H_0\diff(S_{\bullet},A)\rw 0
\end{split}
\end{equation}
\begin{equation}\label{class.eq17}
\begin{split}
   & 0\rw H^0\der(S_{\bullet},A)\rw\der(G,A)\rw\der(\Phi,A)\rw H^1\der(S_{\bullet},A)\rw \\
   & \rw H^2(G,A)\rw D^1(\Phi,A)\rw \cdots.
\end{split}
\end{equation}  \enpr

\medskip

Since, from Proposition \ref{ascro.pro1}, $D_n(\Phi,\zbb)\cong\zeta H^{CCG}_{n+1}\tgmu$,
the long exact homology sequences of Theorem \ref{class.the1b} for the case $A=\zbb$
recovers the result of \cite[Corollary 4]{glp} which is established there via a different
method . The following are consequences of the theorem above.
\begin{corollary}\label{class.cor1}
    Let $\Phi=(T,G,\mu)$ and $\Phi'=(T',G',\mu')$ be two crossed modules acting on the
    abelian group $A$. Suppose that there is a weak equivalence $\Phi\rw \Phi'$ inducing isomorphisms $H_*(G,A)\cong
    H_*(G',A)$, $\;H^*(G,A)\cong H^*(G',A)$. Then for each $n\geq 0$
\begin{equation*}
  D_n(\Phi,A)\cong D_n(\Phi',A),\qquad D^n(\Phi,A)\cong D^n(\Phi',A).
\end{equation*}
\end{corollary}
\prf By hypothesis there exists a crossed module homomorphism $(f,g):\Phi\rw \Phi'$
inducing isomorphisms of homotopy groups. Let
$\gbb_{\bullet}\Phi=(T_{\bullet},G_{\bullet},\mu_{\bullet})$,
$\;\gbb_{\bullet}\Phi'=(T'_{\bullet},G'_{\bullet},\mu'_{\bullet})$,
$\;S_{\bullet}=G_{\bullet}/\mu_{\bullet}(T_{\bullet})$,
$\;S'_{\bullet}=G'_{\bullet}/\mu'_{\bullet}(T'_{\bullet})$. The homomorphisms
$\gbb_{\bullet}(f,g):\gbb_{\bullet}\Phi\rw\gbb_{\bullet}\Phi'$ induce homomorphisms
$S_{\bullet}\rw S'_{\bullet}$. By Lemma \ref{ascro.lem1} we have a commutative diagram of
chain complexes
\begin{diagram}[h=1.5em,s=2em]
    0 & \rTo & A\otund{\zbb\p}\diff\gbb_{\bullet}\Phi & \rTo & A\otund{\zbb\p}\diff
    G_{\bullet}& \rTo & A\otund{\zbb\p}\diff S_{\bullet} & \rTo & 0 \\
    && \dTo && \dTo && \dTo \\
    0 & \rTo & A\otund{\zbb\pp}\diff\gbb_{\bullet}\Phi' & \rTo & A\otund{\zbb\pp}\diff
    G'_{\bullet}& \rTo & A\otund{\zbb\pp}\diff S'_{\bullet} & \rTo & 0\,,
\end{diagram}
where $\p$ and $\pp$ are the first homotopy groups of $\Phi$ and $\Phi'$ respectively.
Since $\Phi$ and $\Phi'$ are weakly equivalent, by Theorem \ref{class.the1b}
$H_*(A\otund{\zbb\p}\diff S_{\bullet})\cong H_*(A\otund{\zbb\pp}\diff S'_{\bullet})$.
Taking the corresponding long exact homology sequence in each row of the above diagram
and using the hypothesis that $H_*(G,A)\cong H_*(G',A)$ we deduce by the Five Lemma that
$D_*(\Phi,A)\cong D_*(\Phi',A)$. The case of cohomology is similar.\enpr
\medskip
\begin{corollary}\label{class.cor2}
    Let $\Phi=(T,G,\mu)$ be a crossed module acting on the abelian group $A$, and
    $\Phi'=(R,F,\delta)$ a crossed module weakly equivalent to $\Phi$ with
    $F$ a free group. Let $\p\cong G/\mu(T)\cong F/\delta(R)$. Then
\begin{equation*}
  H_{n+2}(B(\Phi),A)\cong D_n(\Phi',A),\qquad H^{n+2}(B(\Phi),A)\cong D^n(\Phi',A)
\end{equation*}
for each $n\geq{1}$ and
\begin{align*}
  H_2(B(\Phi),A)\cong & \ker(A\otund{\zbb \p}\diff \Phi'\rw A\otund{\zbb \p}\diff
  F)\\
    H^2(B(\Phi),A)\cong & \mathrm{coker}\:(\der(F,A)\rw
    \der(\Phi',A)).
\end{align*}
\end{corollary}
\smallskip
\prf $\;$Apply Theorem \ref{class.the1b} using the fact that $\;H_*(B(\Phi),A)\;\cong\;
H_*(B(\Phi'),A)$, $\;H^*(B(\Phi),A)\cong H^*(B(\Phi'),A)$ and $H_n(F,A)=H^n(F,A)=0$ for
$n\geq 2$.\enpr

\medskip

Notice that, since
$\der(\Phi,A)\cong\{D\in\hund{\mathcal{G}r}(T,A)\;|\;^{g}D(t)=D(\,^{g}t)\}$ from
Corollary \ref{class.cor2} we recover the result of \cite[Theorem 6]{elli}, which is
established there by different method. Moreover, taking $A=\zbb$ and using the fact that
$\zbb\otund{\zbb\p}\diff\Phi'\cong R/[R,F]$ and $\,\zbb\otund{\zbb\p}\diff F\,\cong\,
F/[F,F]$, we recover the Hopf-type formula for $H_2(B(\Phi),\zbb)$ first proved in
\cite[Theorem 6]{elli}.

\medskip

Our aim in the remaining part of this section is to give a description of the
(co)homology $D_*(\Phi,A)$ and $D^*(\Phi,A)$ without using cotriples.

Let $\Phi=\tgmu$ be a crossed module. The inclusion $(1,G,i)\hookrightarrow\tgmu$ induces
an inclusion of simplicial groups $N_*^{-1}(1,G,i)\rw N_*^{-1}\tgmu$. In turn this
determines for each $m,n$ an injection $\;C_m(N_n^{-1}(1,G,i),A)\remb
C_m(N_n^{-1}\tgmu,A)$ and a surjection $\;C^m(N_n^{-1}\tgmu,A)\ronto
C^m(N_n^{-1}(1,G,i),A)$. We therefore have an injective map of chain complexes of abelian
groups
\begin{equation}\label{class.eq1}
    \tot C_*(N_*^{-1}(1,G,i),A)\remb\tot C_*(N_*^{-1}\tgmu,A)
\end{equation}
and a surjective map of cochain complexes of abelian groups
\begin{equation}\label{class.eq2}
    \tot C^*(N_*^{-1}\tgmu,A)\ronto \tot C^*(N_*^{-1}(1,G,i),A).
\end{equation}
Denote by $\beta_{\bullet}(\Phi,A)$ the cokernel in (\ref{class.eq1}) and by
$\beta^{\bullet}(\Phi,A)$ the kernel in (\ref{class.eq2}).

\begin{lemma}\label{class.lem1}
    Let $\Phi=\tgmu$ be an aspherical crossed module with $T,G,\p$ free
    groups, and let $A$ be a $\p$-module. Then
\begin{equation*}
\begin{split}
   & H_n\beta_\bullet(\Phi,A)=H^n\beta^\bullet(\Phi,A)=0\quad\text{for each}\quad n\neq 2, \\
   & H_2\beta_\bullet(\Phi,A)\cong\diff(\Phi,A),\\
   & H^2\beta^\bullet(\Phi,A)\cong\der(\Phi,A).
\end{split}
\end{equation*}
\end{lemma}

\prf Consider the long exact homology sequence associated to the short exact sequence of
chain complexes
\begin{equation*}
\tot C_*(N_*^{-1}(1,G,i),A)\remb\tot C_*(N_*^{-1}\tgmu,A) \ronto\beta_\bullet(\Phi,A).
\end{equation*}
We have $H_n(G,A)=0$ for $n\geq 2$ as $G$ is free; since $\tgmu$ is weakly equivalent to
$(1,\p,i)$ and $\p$ is free, $H_n(B(\Phi),A)\cong H_n(\p,A)=0$ for $n\geq 2$ and
$\!H_1(B(\Phi),A)\!\cong \!H_1(\p,A)$. $\!$Thus this long exact homology sequence gives
$\!H_n\beta_\bullet(\!\Phi,\!A)\!=\!H^n\beta^\bullet(\!\Phi,\!A)\!=0$ for $n>2$ and the
exact sequence
\begin{equation*}
\begin{split}
   & 0\rw H_2\beta_\bullet(\Phi,A)\rw H_1(G,A)\rw H_1(\p,A)\rw H_1\beta_\bullet(\Phi,A)\rw \\
   & \rw H_0(G,A)\rw H_0(\p,A)\rw H_0\beta_\bullet(\Phi,A)\rw 0.
\end{split}
\end{equation*}
Since $H_0(G,A)\cong A_G\cong A_{\p}\cong H_0(\p,A)$ it follows that
$H_i\beta_\bullet(\Phi,A)=0$ for $i=0,1$ and
\begin{equation*}
  H_2\beta_\bullet(\Phi,A)\cong\ker(H_1(G,A)\rw H_1(\p,A)).
\end{equation*}
On the other hand, since $\p$ is free, the five-term homology sequence associated to the
extension $T\remb G\ronto\p$ reduces to
\begin{equation*}
  0 \rw A\otund{\zbb\p}T_{ab}\rw H_1(G,A)\rw H_1(\p,A)\rw 0.
\end{equation*}
Thus by Lemma \ref{ascro.lem1}, $\;H_2\beta_\bullet(\Phi,A)\cong\diff(\Phi,A)$. The
argument for cohomology is similar.\enpr
\medskip
\begin{proposition}\label{class.pro1}
    Let $\Phi$ be a crossed module acting on the abelian group $A$. Then for each $n\geq2$
\begin{align*}
   & H_n(\beta_{\bullet}(\Phi,A))\cong D_{n-2}(\Phi,A),\qquad
   H^n(\beta^{\bullet}(\Phi,A))\cong D^{n-2}(\Phi,A).
\end{align*}
\end{proposition}
\medskip
\prf Let $\gbb_{\bullet}\Phi=(T_{\bullet},G_{\bullet},i_{\bullet})$ and consider the
bicomplexes $\{\psi_{pq}\}$, $\{\chi_{pq}\}$, $\{\mathcal{L}_{pq}\}$,
\begin{align*}
  \psi_{pq} & =\beta_q((T_p,G_p,i_p),A), \\
  \chi_{pq} & =(\tot C_*(N^{-1}_*(T_p,G_p,i_p),A))_q,\\
  \mathcal{L}_{pq} & = (\tot C_*(N^{-1}_*(1,G_p,i),A))_q.
\end{align*}
We aim to show that for each $n\geq 2$
\begin{equation}\label{class.eq2b}
  H_n\tot\psi_{\bullet\bullet}\cong H_n\beta_{\bullet}(\Phi,A).
\end{equation}
The morphism of simplicial crossed modules
$(T_{\bullet},G_{\bullet},i_{\bullet})\!\rw\!\tgmu$ and $(1,G_{\bullet},i)\!\rw\!(1,G,i)$
(here $\tgmu$ and $(1,G,i)$ are thought of as constant simplicial crossed modules) induce
morphisms of double complexes $\chi_{\bullet\bullet}\rw\tot C_*(N^{-1}_*\tgmu,A)$ and
$\mathcal{L}_{\bullet\bullet}\rw\tot C_*(N^{-1}_*(1,G,i),A)$. We claim that these
morphisms induce isomorphisms in the total homologies in dimensions $n\geq 2$. In fact if
$S_{\bullet}=G_{\bullet}/i_{\bullet}(T_{\bullet})$, the double complex
$\chi_{\bullet\bullet}$ gives rise to a spectral sequence
\begin{align*}
   & E^1_{pq}=H_q\tot C_*(N^{-1}_*(T_p,G_p,i_p),A)\cong H_q(B(T_p,G_p,i_p),A)\cong \\
   & \cong H_q(B(1,S_p,i),A)\cong H_q(S_p,A),\\
   & E^2_{pq}=H_pE^1_{*q}.
\end{align*}
On the other hand since each $S_p$ is a free group, the double complex
 $\alpha_{\bullet\bullet}$
\begin{diagram}[h=1.5em,s=1.6em]
    && 0 && 0 && 0 && 0 &&\\
    && \dTo && \dTo && \dTo && \dTo && \\
    \cdots&\rTo & A\otund{\zbb S_3}\ifr_{S_3} & \rTo & A\otund{\zbb S_2}\ifr_{S_2} &
    \rTo & A\otund{\zbb S_1}\ifr_{S_1} & \rTo & A\otund{\zbb S_0}\ifr_{S_0} & \rTo & 0\\
    && \dTo && \dTo && \dTo && \dTo\\
    \cdots&\rTo & A & \rTo_{\id} & A & \rTo_{\id} & A & \rTo_{\id} & A & \rTo & 0 \\
    && \dTo && \dTo && \dTo && \dTo && \\
    && 0 && 0 && 0  && 0 &&.
\end{diagram}
gives rise to a spectral sequence with the same $E^1$ and $E^2$ terms, hence
$H_n\tot\alpha_{\bullet\bullet}\cong H_n\tot\chi_{\bullet\bullet}$. As shown in the proof
of theorem \ref{class.the1b}, for each $n\geq 2$
$H_n\tot\alpha_{\bullet\bullet}=H_{n-1}(A\otund{\zbb S_{\bullet}}\ifr_{S_{\bullet}})$ and
by Theorem \ref{class.the1b} ii) $H_{n-1}(A\otund{\zbb
S_{\bullet}}\ifr_{S_{\bullet}})\cong H_n(B\Phi,A)$; hence
$H_n\tot\chi_{\bullet\bullet}\cong H_n(B\Phi,A)\cong H_n\tot C_*(N^{-1}_*\tgmu,A)$ for
each $n\geq 2$.

Similarly, since each $G_p$ is free,
\begin{equation*}
  \tot C_*(N^{-1}_*(1,G_{\bullet},i),A)_q=
  \begin{cases}
    A & q=0, \\
    A\oplus(A\otund{\zbb G_{\bullet}}\ifr_{G_{\bullet}}) & q>0.
  \end{cases}
\end{equation*}
Therefore $\mathcal{L}_{\bullet\bullet}$ gives rise to a spectral sequence
\begin{equation*}
  E^1_{pq}=H_p(\tot C_*(N^{-1}_*(1,G_{\bullet},i),A))_q=
  \begin{cases}
    0 & q=0, \\
    H_p(A\otund{\zbb G_{\bullet}}\ifr_{G_{\bullet}}) & q>0.
  \end{cases}
\end{equation*}
\begin{equation*}
  E^2_{pq}=H_q E^1_{p*}=
  \begin{cases}
    0 & q\neq 1, \\
    H_p(A\otund{\zbb G_{\bullet}}\ifr_{G_{\bullet}}) & q=1.
  \end{cases}
\end{equation*}
It follows that, for each $n\geq 2$,
$\;H_n\tot\mathcal{L}_{\bullet\bullet}=E^2_{n-1,1}=H_{n-1}(A\otund{\zbb
G_{\bullet}}\ifr_{G_{\bullet}})=H_n(G,A)\cong H_n\tot C_*(N^{-1}_*(1,G,i),A)$. This
proves the claim.

Consider the commutative diagram of short exact sequences of double complexes of abelian
groups
\begin{diagram}[h=1.8em]
    \mathcal{L}_{\bullet\bullet} & \rEmbed & \chi_{\bullet\bullet} & \rOnto &
    \psi_{\bullet\bullet}\\
    \dTo && \dTo && \dTo\\
    \tot C_*(N^{-1}_*(1,G,i),A) & \rEmbed & \tot C_*(N^{-1}_*\tgmu,A) & \rOnto &
    \beta_{\bullet}(\tgmu,A) .
\end{diagram}
Taking the induced long exact sequences in total homologies in each row of the diagram,
from the claim and the Five lemma it follows that, for each $n\geq 2$
\begin{equation*}
  H_n\tot\psi_{\bullet\bullet}\cong H_n\beta_{\bullet}(\tgmu,A)
\end{equation*}
which is (\ref{class.eq2b}).

On the other hand, by Lemma \ref{class.lem1} the double complex $\psi_{\bullet\bullet}$
gives rise to a spectral sequence
\begin{equation*}
  E^1_{pq}=H_q\beta_{\bullet}((T_p,G_p,i_p),A)=
  \begin{cases}
    0 & q\neq{2}, \\
    \diff((T_p,G_p,i_p),A) & q=2
  \end{cases}
\end{equation*}
so that
\begin{equation*}
  E^2_{pq}=H_pE^1_{*q}=
  \begin{cases}
    0 & q\neq 2, \\
    H_p\diff((T_{\bullet},G_{\bullet},i_{\bullet}),A)=D_p(\Phi,A) & q=2.
  \end{cases}
\end{equation*}
Hence $E^2_{pq}\Rightarrow
H_{p+q}\mathrm{Tot}\psi_{{\bullet}{\bullet}}=H_{p+q}\beta_{\bullet}((T,G,\mu),A)$
collapses, giving
\begin{equation*}
  H_n\beta_{\bullet}(\Phi,A)=E^2_{n-2,2}=D_{n-2}(\Phi,A)
\end{equation*}
for $n\geq 2$. The argument for cohomology is similar.\enpr

\medskip

We notice that the description of the (co)homology $D^*(\Phi,A)$ and $D_*(\Phi,A)$ given
in the above proposition gives rise to a version of Theorem \ref{class.the1b} (ii), which
differs from the previous one in low dimensions.
\begin{corollary}\label{class.cor1b}
    Let $\Phi=\tgmu$ be a crossed module, $A$ a $\p$-module. There exist long exact
    (co)homology sequences
       \begin{equation*}
  \begin{split}
  & \rw D_n(\Phi,A)\rw H_{n+1}(G,A)\rw H_{n+1}(B(T,G,\mu),A)\rw D_{n-1}(\Phi,A)\rw\\
  &  \cdots \rw H_2(G,A)\rw H_2(B(T,G,\mu),A)\rw A\otund{\zbb \p}\diff\Phi\rw\\
  & \rw H_1(G,A)\rw H_1(\p,A)\rw 0
  \end{split}
\end{equation*}
\begin{equation*}
  \begin{split}
  & 0\rw H^1(\p,A)\rw H^1(G,A)\rw \der(\Phi,A)\rw H^2(B(T,G,\mu),A)\rw \\
  & \rw H^2(G,A)\rw D^1(\Phi,A)\rw \cdots
  \end{split}
\end{equation*}
\end{corollary}
\prf Take the long exact (co)homology sequences associated to the short exact sequences
of (co)chain complexes
\begin{equation*}
\begin{split}
   & \mathrm{Tot\,}C_*(N_*^{-1}(1,G,i),A) \remb \mathrm{Tot\,}C_*(N_*^{-1}
   \tgmu,A)\ronto\beta_{\bullet}(\Phi,A) \\
   & \beta^{\bullet}(\Phi,A)\remb\mathrm{Tot\,}C^*(N_*^{-1}\tgmu,A)\ronto\mathrm{Tot}\,
   C^*(N_*^{-1}(1,G,i),A).
\end{split}
\end{equation*}
Apply Proposition \ref{class.pro1} and the fact (see \cite{elli}) that
$H_1(B(\Phi),A)\cong H_1(\p,A)$ and $H^1(B(\Phi),A)=H^1(\p,A)$.\enpr

\medskip

We finally notice that a more topological approach than the one given in this paper
should allow to obtain a topological interpretation of the (co)homology $D^*(\Phi,A)$ and
$D_*(\Phi,A)$ as relative (co)homology of the pair of spaces $(B\Phi,BG)$ with local
coefficients. In fact, it is reasonable to conjecture that, for each $n>0$,
$\;H^n\beta^\bullet(\Phi,A)$ (resp. $H_n\beta_{\bullet}(\Phi,A)$) is isomorphic to
$H^n(B\Phi,BG;A)$ (resp. $H_n(B\Phi,$ $BG;A)$), so that by Proposition \ref{class.pro1},
for each $n\geq 2$ $D^{n-2}(\Phi,A)$ (resp. $D_{n-2}(\Phi,A)$) would be isomorphic to
$H^n(B\Phi,BG;A)$ (resp. $H_n(B\Phi,BG;A)$).

%%%%%% New section %%%%%%%%%%%%%%%%%%%%%%%%%%%%%%%%%%%%%%%%%%%%%%%%%%%%%%%%%%%%%%%%%%%
\section{An example}\label{exa}
Let $M$ be a \md{G} and consider the crossed module $\Phi=(M,G,0)$. The map of crossed
modules $(i,\id_G)\!:(1,G,i)\rw (M,G,0)$ has a section $(0,\id_G):(M,G,0)\rw(1,G,i)$.
Therefore the corresponding map $\,B(G)\,\hookrightarrow\, B(M,G,0)\, $ has a section
$\,B(M,G,0)\rw B(G)$. Hence the long exact sequences of Theorem \ref{class.the1b} give
split short exact sequences for each $n\geq 2$
\begin{equation}\label{exa.eq1a}
\begin{split}
  & 0 \rw H_n(G,A)\;\pile{\leftarrow \\ \rw}\; H_n(B(M,G,0),A)\rw D_{n-2}((M,G,0),A)\rw
  0\\
  & 0 \rw D^{n-2}((M,G,0),A)\rw H^n (B(M,G,0),A)\;\pile{\leftarrow \\ \rw}\; H^n(G,A)\rw
  0.
\end{split}
\end{equation}
It follows that for each $n\geq 2$
\begin{equation}\label{exa.eq1}
\begin{split}
   & H_n(B(M,G,0),A)\cong H_n(G,A)\oplus D_{n-2}((M,G,0),A) \\
   & H^n(B(M,G,0),A)\cong H^n(G,A)\oplus D^{n-2}((M,G,0),A).
\end{split}
\end{equation}
Recall that for every crossed module $\tgmu$ with homotopy groups $\p$ and $\pi_2$ there
is a fibration sequence
\begin{equation*}
  K(\pi_2,2)\rw |B\tgmu |\rw K(\pi_1,1)
\end{equation*}
where $K(\pi_2,2)$ and $K(\p,1)$ are Eilenberg--MacLane spaces. Hence we have
corresponding Serre spectral sequences:
\begin{equation*}
\begin{split}
   & E^2_{pq}=H_p(\pi_1,H_q(K(\pi_2,2),A))\Rightarrow H_{p+q}(B\tgmu,A) \\
   & E_2^{pq}=H^p(\pi_1,H^q(K(\pi_2,2),A))\Rightarrow H^{p+q}(B\tgmu,A).
\end{split}
\end{equation*}
In the following proposition we shall use the well known fact that, for every abelian
group $A$
\begin{equation}\label{exa.eq2}
  H_1 K(A,2)=0=H_3 K(A,2),\qquad H_2 K(A,2)\cong A,\qquad H_4 K(A,2)\cong\Gamma^2A
\end{equation}
where $\Gamma^2$ denotes Whitehead's universal quadratic functor.
\begin{proposition}\label{exa.pro1}
    Let $M$ be a \md{\zbb G} and let $\Phi=(M,G,0)$ act on the abelian group $A$.
    \\a) There are exact sequences
\begin{equation*}
\begin{split}
   & D_2(\Phi,A)\rw H_2(G,M\otmz A)\rw H_0(G,\tor^{\zbb}_{1}(M,A))\rw \qquad\qquad\;\;\\
   & \rw D_1(\Phi,A)\rw H_1(G,M\otmz A)\rw 0\\
\end{split}
\end{equation*}
\begin{equation*}
\begin{split}
   & 0\rw H^1(G,\homz(M,A))\rw D^1(\Phi,A)\rw H^0(G,\ext^1_\zbb(M,A))\rw\\
   & \rw H^2(G,\homz(M,A))\rw D^2(\Phi,A).
   \end{split}
\end{equation*}
b) If $\,\tor^\zbb_1(M,A)=0$ then there is an exact sequence
\begin{equation*}
\begin{split}
   & D_3(\Phi,A)\rw H_3(G,M\otmz A)\rw H_0(G,\Gamma^2 M\otmz A)\rw \\
   & \rw D_2(\Phi,A)\rw H_2(G,M\otmz A)\rw 0.
\end{split}
\end{equation*}
$\quad $If $\,\ext^1_\zbb(M,A)=0$ then there is an exact sequence
\begin{equation*}
\begin{split}
   & 0\rw H^2(G,\homz(M,A))\rw D^2(\Phi,A)\rw H^0(G,\homz(\Gamma^2 M,A))\rw \\
   & \rw H^3(G,\homz(M,A))\rw D^3(\Phi,A).
\end{split}
\end{equation*}
\end{proposition}
\prf The Serre spectral sequences for the crossed module $(M,G,0)$ are
\begin{align*}
   & E^2_{pq}=H_p(G,H_q(K(M,2),A))\Rightarrow H_{p+q}(B(M,G,0),A) \\
   & E_2^{pq}=H^q(G,H^q(K(M,2),A))\Rightarrow H^{p+q}(B(M,G,0),A)
\end{align*}
while the same spectral sequence for the crossed module $(1,G,i)$ has
$E^2_{pq}=E^{pq}_2=0$ for $q\neq 0$; also $H_n(B(1,G,i),A)=H_n(G,A)$ and
$H^n(B(1,G,i),A)=H^n(G,A)$.

Hence from (\ref{exa.eq1a}) we deduce that there are spectral sequences
\begin{equation*}
\begin{split}
  &\widetilde{E}^2_{pq}=H_p(G,H_{q+2}(K(M,2),A))\Rightarrow D_{p+q}((M,G,0),A)\\
  &\widetilde{E}_2^{pq}=H^p(G,H^{q+2}(K(M,2),A))\Rightarrow D^{p+q}((M,G,0),A).
\end{split}
\end{equation*}
The corresponding exact sequences of low degree terms are
\begin{equation*}
\begin{split}
   & D_2(\Phi,A)\rw H_2(G,H_2(K(M,2),A))\rw H_0(G,H_3(K(M,2),A))\rw \qquad\\
   & \rw D_1(\Phi,A)\rw H_1(G,H_2(K(M,2),A))\rw 0.
\end{split}
\end{equation*}
\begin{equation*}
\begin{split}
   & 0\rw H^1(G,H^2(K(M,2),A))\rw D^1(\Phi,A)\rw H^0(G,H^3(K(M,2),A))\rw \\
   & \rw H^2(G,H^2(K(M,2),A))\rw D^2(\Phi,A).
\end{split}
\end{equation*}
By (\ref{exa.eq2}) and the universal coefficient theorem,
\begin{equation*}
\begin{array}{ll}
  H_2(K(M,2),A)\cong M\otund{\zbb}A,&\quad  H_3(K(M,2),A)\cong\tor^\zbb_1(M,A) \\
  H^2(K(M,2),A)\cong\hund{\zbb}(M,A),&\quad  H^3(K(M,2),A)\cong\ext^1_\zbb(M,A)
\end{array}
\end{equation*}
so that part a) follows.

\medskip

If $\tor^\zbb_1(M,A)=0$ then $\widetilde{E}^2_{p1}=H_p(G,\tor^\zbb_1(M,A))=0$. Hence (see
\cite{caei}) there is an exact sequence
$D_3\rw\widetilde{E}^2_{3,0}\rw\widetilde{E}^2_{0,2}\rw D_2\rw\widetilde{E}^2_{2,0}\rw
0$, i.e.
\begin{equation*}
\begin{split}
   & D_3(\Phi,A)\rw H_3(G,H_2(K(M,2),A)\rw H_0(G,H_4(K(M,2),A))\rw \\
   & \rw D_2(\Phi,A)\rw H_2(G,H_2(K(M,2),A))\rw 0.
\end{split}
\end{equation*}
From (\ref{exa.eq2}) and the universal coefficient theorem, we have
$H_4(\!K\!(M,2)\!,\!A\!)=\Gamma^2M\otund{\zbb}A$ so that part b) follows for the homology
case.

\medskip

If  $\,\ext^1_\zbb(M,A)=0$ then $\widetilde{E}^{p1}_2=H^p(G,\ext^1_\zbb (M,A))=0$ so
there is an exact sequence $0\rw \widetilde{E}^{2,0}_2\rw
D^2\rw\widetilde{E}_2^{0,2}\rw\widetilde{E}_2^{3,0}\rw D^3\,$,  i.e.
\begin{equation*}
\begin{split}
   & 0\rw H^2(G,H^2(K(M,2),A))\rw D^2(\Phi,A)\rw H^0(G,H^4(K(M,2),A))\rw \\
   & \rw H^3(G,H^2(K(M,2),A))\rw D^3(\Phi,A).
\end{split}
\end{equation*}
From the universal coefficient theorem $H^4(K(M,2),A))=\homz(\Gamma^2 M,A)$ so that part
b) follows.\enpr
\medskip

We point out that for the case of homology with $\zbb$-coefficients, (\ref{exa.eq1}) and
the homology exact sequences of Proposition \ref{exa.pro1} were also given in \cite{glp}.

%%%%%% New section %%%%%%%%%%%%%%%%%%%%%%%%%%%%%%%%%%%%%%%%%%%%%%%%%%%%%%%%%%%%%%%%%%%
\section{An application to the cohomology of the classifying space}\label{app}
In this section we apply the simplicial description of the cohomology of the classifying
space of a crossed module proved in Theorem \ref{class.the1b} ii) to give an
interpretation of these cohomology groups in dimensions $n=2,3$.
\begin{lemma}\label{app.lem1}
    Let $\Phi=\tgmu$ be a crossed module, $A$ an abelian group. Then $(1,A,i)$ is a
    \md{\Phi} if and only if $A$ is a \md{\p}. In this case there is an isomorphism
\begin{equation*}
  \der(\Phi,(1,A,i))\cong\der(\p,A).
\end{equation*}
\end{lemma}
\prf From Section 1, the singular object $(1,A,i)$ is a \md{\Phi} if and only if there is
a split extension of crossed modules
\begin{equation}\label{app.eq1}
  (1,A,i)\remb(T', G',\mu')\;\pile{\ronto \\ \leftarrow}
  \;\tgmu.
\end{equation}
In particular we have split short exact sequences of groups $1\remb T'\ronto T$ and
$A\remb G'\;\pile{\ronto\\ \leftarrow}\; G$; thus we can assume that $T'=T$, $\;G'\cong
A\rt G$, that the map $\pr_G:A\rt G\rw G$ is $\pr_G(a,g)=g$ and that the map $i_G:G\rw
A\rt G$ is $i_G(g)=(0,g)$. Hence we have the split short exact sequence
\begin{equation}\label{app.eq1a}
  (1,A,i)\rEmbed (T,A\rt G,\mu')\pile{\lTo^{(\id_T,i_G)}\\ \rOnto_{(\id_T,\pr_G)}} \tgmu.
\end{equation}
Since $(\id_T,i_G)$ is a crossed module map, $\mu'\id_T=i_G\mu$, hence
$\mu'(t)=(0,\mu(t))$ for all $t\in T$, that is $\mu'=(0,\mu)$. Since $(\id_T,\pr_G)$ is a
crossed module map, $^{(a,g)}t=\:^{g}t$ for all $(a,g)\in A\rt G, \;\;t\in T$. The axioms
of crossed module for $(T,A\rt G, (0,\mu))$ give, for all $(a,g)\in A\rt G,\;\;t\in T$
\begin{equation*}
  (0,\mu)(\:^{(a,g)}t)=(a,g)(0,\mu(t))(a,g)^{-1}.
\end{equation*}
An easy calculation shows that this is equivalent to
\begin{equation*}
  (0,\mu(\:^{g}t))=(a-\:^{\mu(\:^{g}t)}a,\mu(\:^{g}t))
\end{equation*}
for all $t\in T$, $g\in G$, $a\in A$. It follows that $a=\:^{\mu(\:^{g}t)}a$. In
particular, taking $g=1$ we obtain $a=\:^{\mu(t)}a$ for all $t\in T$, $a\in A$, so that
$A$ is a $\p$-module.

Conversely, if $A$ is a \md{\p}, then (\ref{app.eq1a}) is a split singular extension of
crossed modules, hence $(1,A,i)$ is a \md{\Phi}. We have
\begin{equation*}
  \der(\Phi,(1,A,i))\cong\hund{\cm/\Phi}(\Phi,(T,A\rtimes G,(0,\mu)).
\end{equation*}
We now show that there is an isomorphism
\begin{equation*}
  \alpha:\hund{\cm/\Phi}(\Phi,(T,A\rtimes G,(0,\mu))\rw\der(\p,A).
\end{equation*}
Let $\alpha(\id_T,(D,\id_G))=D$. Then $D\in\der(G,A)$ and since $(D,\id_G)\mu=(0,\mu)$ we
have $D\mu=0$ so that $D\in\der(\p,A)$. Clearly $\alpha$ is injective. Given
$D\in\der(\p,A)$ let $\overline{D}(g)=D(g\mu(T)),\; g\in G$. Then
$\overline{D}\in\der(G,A)$ and $\alpha(\id_T,( \overline{D},\id_G))=D$, so that $\alpha$
is also surjective.\enpr
\medskip

We can similarly define a functor $\der(\mi,(1,A,i)):\cm/\Phi\rw\ab$ on the slice
category.
\medskip

In the next proposition we show that the cohomology of the classifying space of a crossed
module can be described as cotriple cohomology.
\begin{proposition}\label{app.pro1}
    Let $\Phi=\tgmu$ be a crossed module, $A$ a \md{\p}. Then for each $n>0$
\begin{equation*}
  H^n(B(\Phi),A)\cong H^{n-1}\der(\gbb_{\bullet}\Phi,(1,A,i)).
\end{equation*}
\end{proposition}
\prf Let $\gbb_{\bullet}\Phi=(T_{\bullet},G_{\bullet},i_{\bullet})$,
$\;S_{\bullet}=G_{\bullet}/i_{\bullet}(T_{\bullet})$. From Lemma \ref{app.lem1}
\begin{equation*}
  H^n\der(\gbb_{\bullet}\Phi,(1,A,i))\cong H^n\der(S_{\bullet},A).
\end{equation*}
The result follows from Theorem \ref{class.the1b} ii).\enpr

\medskip

The following corollary generalizes a result of \cite[Theorem 10 (iv)]{ccg} which is
established there in the case of aspherical crossed modules.
\begin{corollary}\label{app.cor1}
    Let $\Phi=\tgmu$ be a crossed module, $A$ a trivial \md{\p}. Then for each $n>0$
\begin{equation*}
  H^n_{CCG}(\tgmu,(1,A,i))\cong H^n(B(\Phi),A).
\end{equation*}
\end{corollary}
\prf Let $\gbb_{\bullet}\Phi=(T_{\bullet},G_{\bullet},i_{\bullet})$,
$\;S_{\bullet}=G_{\bullet}/i_{\bullet}(T_{\bullet})$. Since actions are trivial
\begin{equation*}
  \der(S_{\bullet},A)\cong\hund{\mathcal{G}p}(S_{\bullet},A)\cong
  \hund{\cm}((T_{\bullet},G_{\bullet},i_{\bullet}),(1,A,i)).
\end{equation*}
By definition $H^n_{CCG}(\tgmu,(1,A,i))\cong
H^{n-1}\hund{\cm}((T_{\bullet},G_{\bullet},i_{\bullet}),(1,A,i))$ and the result follows
from Theorem \ref{class.the1b} ii).\enpr
\medskip

We finally obtain the interpretation for the second and third cohomology group of the
classifying space. We need the notion of singular and 2-fold special extensions of
$\tgmu$ by $(1,A,i)$.
\begin{definition}\label{app.def1b}
    Let $\Phi\!=\!\tgmu\!$ be a crossed module acting on an abelian group $A$.
\begin{description}
  \item[i)] A singular extension of $\,\tgmu\,$ by $\,(1,A,i)\,$ is a short exact sequence of
  crossed modules
\begin{equation}\label{app.eq1b}
  (1,A,i)\remb(T,G',\mu')\overset{(\id_T,f)}{\ronto}\tgmu
\end{equation}
such that the corresponding short exact sequence of cat$^1$-groups
\begin{equation}\label{app.eq2b}
  (1\times A,\id,\id)\remb(T\rt G',d',s')\ronto(T\rt G,d,s)
\end{equation}
is a singular extension of $(T\rt G,d,s)$ by the $(T\rt G,d,s)$-module $(1\times
A,\id,\id)$ in the sense of categories of interest \cite{vale}.
  \item[ii)] A 2-fold special extension of $\,\tgmu\,$ by $\,(1,A,i)\,$ is an exact sequence
  of crossed modules
\begin{equation}\label{app.eq3b}
  (1,A,i)\overset{i}{\remb}(T'',G'',\mu'')\xrw{(\alpha,\beta)}(T',G',\mu')
  \overset{(f,r)}{\ronto}\tgmu
\end{equation}
such that the corresponding exact sequence of cat$^1$-groups
\begin{equation}\label{app.eq4b}
  (1\times A,\id,\id)\remb(T''\rt G'',d'',s'')\xrw{(\alpha,\beta)}(T'\rt
  G',d',s')\overset{(f,r)}{\ronto}(T\rt G,d,s)
\end{equation}
is a 2-fold special extension of $(T\rt G,d,s)$ by the $(T\rt G,d,s)$-module $(1\times
A,\id,\id)$ in the sense of categories of interest \cite{vale}.
\end{description}
\end{definition}
A more explicit characterization of singular and 2-fold special extensions of $\tgmu$ by
$(1,A,i)$ can be given as follows.
\begin{lemma}\label{app.lem1b}
    Let $\Phi=\tgmu$ be a crossed module acting on the abelian group $A$.
\begin{description}
  \item[i)] A singular extension of $\tgmu$ by $(1,A,i)$ consists of a short exact
  sequence of crossed modules (\ref{app.eq1b}) such that if $f':G\rw G'$ is a set map with
  $ff'=\id_G$, it is $f'(g)af'(g^{-1})=\:^{[g]}a$ for all $g\in G$, $a\in A$ where $[g]=g\mu(T)\in
  \p$ and $^{[g]}a$ is the given $\p$-module action on $A$.
  \item[ii)] A 2-fold special extension of $\tgmu$ by $(1,A,i)$ consists of an
  exact sequence of crossed modules (\ref{app.eq3b}) where
\begin{equation}\label{app.eq4c}
\begin{diagram}[s=1.8em]
    T'' & \rTo^\alpha & T' &&&& h:T'\times G''\rw T''\\
    \dTo^\mu && \dTo_{\mu'} \\
    G'' & \rTo_\beta & G'
\end{diagram}
\end{equation}
is a crossed square and if $r':G\rw G'$ is a set map with $rr'=\id_G$, then for all $g\in
G$ , $a\in A$,
\begin{equation*}
  ^{r'(g)}a=\:^{[g]}a.
\end{equation*}
Here $[g]=g\mu(T)\in \p$, $\:^{[g]}a$ is the given $\p$-module action on $A$ and
$^{r'(g)}a$ is the action of $G'$ on $G''$ in the crossed square (\ref{app.eq4c}).
\end{description}
\end{lemma}
\prf
\begin{description}
  \item[i)] By definition, (\ref{app.eq2b}) is a singular extension of
  cat$^1$-groups. Hence for all $t\in T$, $g\in G$, $a\in A$
\begin{equation*}
  (t,f'(g))(1,a)(t,f'(g^{-1}))=(1, \:^{[g]}a),
\end{equation*}
that is
\begin{equation*}
  (t\:^{f'(g)af'(g^{-1})}t^{-1},f'(g)a f'(g^{-1}))=(1,\:^{[g]}a).
\end{equation*}
Since $(1,A,i)$ is a normal subcrossed module of $\tgmu$,
$t\:^{f'(g)af'(g^{-1})}t^{-1}=1$ for all $t\in T$, $g\in G$, $a\in A$; hence we only
require that $f'(g)a f'(g^{-1})=\:^{[g]}a$.
  \item[ii)] By definition, (\ref{app.eq3b}) is a 2-fold extension in
  cat$^1$-groups. Let $f':T\rw T'$ be a set map with $ff'=\id_T$.
  Since $((T''\rt G'',d'',s''),(T'\rt G',d',s'),(\alpha,\beta))$ is a crossed module
  in the category of cat$^1$-groups, by Lemma \ref{catone.lem1}  (\ref{app.eq4c}) is a crossed
  square and the crossed module action of $T'\rt G'$ on $T''\rt G''$ is given by
\begin{equation}\label{app.eq1c}
  ^{(t',g')}(t'',g'')=(\,^{t'}(\,^{g'}t'')h(t',\,^{g'}g''),\,^{g'}g'').
\end{equation}
Further, arguing as in the proof of Lemma \ref{furpro.lem1b} ii), the induced action of
$(T\rt G,d,s)$ on $(1\times A,\id,\id)$ given by $^{(f'(t),r'(g))}(1,a)$ has to coincide
with the given action, which is $(1,\:^{[g]}a)$. Hence by (\ref{app.eq1c}) we obtain
\begin{equation*}
  (h(f'(t),\:^{r'(g)}a),\:^{r'(g)}a)=(1,\:^{[g]}a)
\end{equation*}
for all $t\in T$, $g\in G$, $a\in A$. By the axioms of crossed squares \cite{lodb}
\begin{equation*}
  \alpha h(f'(t),\:^{r'(g)}a)=f'(t)\;^{\beta(\:^{r'(g)}a)}f'(t)^{-1}=f'(t)f'(t)^{-1}=1.
\end{equation*}
Hence, since $\alpha$ is injective, $h(f'(t),\:^{r'(g)}a)=1$. Therefore we only require
\begin{equation*}
  ^{r'(g)}a=\:^{[g]}a
\end{equation*}
for all $g\in G$, $a\in A$.\enpr
\end{description}

\smallskip

It is possible to introduce an equivalence relation on the set of singular and 2-fold
special extensions of $\tgmu$ by $(1,A,i)$ in a way similar to what explained in \S 4.2.
The sets of equivalence classes of singular and 2-fold special extensions of $\tgmu$ by
$(1,A,i)$ become abelian groups under Baer sum.
\begin{theorem}\label{app.the1}
Let $\,\Phi=\tgmu\,$ be a crossed module, $A$ a \md{\p}. Then $H^2(B(\Phi),A)$ is
isomorphic to the group of equivalence classes of singular extensions of $\tgmu$ by
$(1,A,i)$ and $H^3(B(\Phi),A)$ is isomorphic to the group of equivalence classes of
2-fold special extensions of $\tgmu$ by $(1,A,i)$.
\end{theorem}
\smallskip

\prf $\;\,$From $\;$ Proposition \ref{app.pro1}, $\;H^2\,(\,B(\Phi),\,A)\,\cong\,
H^1\,\der\,(\,\gbb_{\bullet}\,\Phi,\,(1,A,i))$ and $\;H^3(B(\Phi),A)\cong
H^2\der(\gbb_{\bullet}\Phi,(1,A,i))$. By the interpretation in terms of extensions of the
first and second cotriple cohomology groups in categories of interest given in
\cite[Theorems 2.1.3 and 2.2.3]{vale} the result follows.\enpr

%%%%% New Section %%%%%%%%%%%%%%%%%%%%%%%%%%%%%%%%%%%%%%%%%%%%%%%%%%%%%%%%%%
\section{The relationship with cohomology of groups with operators}
Our purpose in this section is to elucidate the relationship between the cohomology
theory $D^*(\tgmu,A)$ of a crossed module $\tgmu$ with coefficients in a $\p$-module $A$
and the cohomology $H^*_G(T,A)$ studied in \cite{cega}. The latter is the cohomology of a
group $T$ endowed with a $G$-action by automorphisms with coefficients in a
$G$-equivariant $T$-module $A$; this consists of an abelian group $A$ with actions of $T$
and $G$ such that
\begin{equation*}
  ^{g}(\,^{t}a)=\,^{^{g}t}(\,^{g}a),\quad g\in G,\;\;t\in T\;\;a\in A.
\end{equation*}
 The possibility that a relationship between the two theories may exist
is suggested by the fact that, by \cite[p. 11]{cega}
\begin{equation*}
  D^0(\tgmu,A)\cong\der(\tgmu,A)\cong \hund{G}(T,A)\cong \der_G(T,A)\cong H^1_G(T,A).
\end{equation*}
We shall also exhibit a counterexample showing that in general $D^n(\tgmu,A)$ and
$H^{n+1}_G(T,A)$ are not isomorphic for $n>0$.

In order to establish the relationship with $H^*_G(T,A)$ we first prove that if $\tgmu$
is a precrossed module, $H^n_G(T,A)$ can be recovered as cohomology of a precrossed
module for $n>0$.

Recall that a \emph{precrossed module} $\tgmu$ consists of a group homomorphism $\mu:T\rw
G$ together with an action of $G$ on $T$ such that $\mu(\,^{g}t)=g\mu(t)g^{-1}$, $\;g\in
G$, $\;t\in T$. A morphism of precrossed modules $(f,h):\tgmu\rw(T',G',\mu')$ consists of
group homomorphisms $f:T\rw T'$, $\;h:G\rw G'$ with $f(\,^{g}t)=\,^{h(g)}f(t)$, $\;t\in
T$, $\;g\in G$. Denote by $\pcm$ the category of precrossed modules. $\pcm$ is equivalent
to the category of \emph{pre-cat$^1$-groups}. A pre-cat$^1$-group is a group $G$ together
with two endomorphisms $d_0,d_1:G\rw G$ such that $d_1d_0=d_0$, $\;d_0d_1=d_1$. A
morphism $f:G\rw G'$ of pre-cat$^1$-groups is a group homomorphism commuting with
$d_0,d_1$.

In \cite{aria} is proved that the category of precrossed modules is tripleable over
$\set$; the corresponding cotriple is then used to define a cotriple (co)homology of
precrossed modules with trivial coefficients. Following the same method used for crossed
modules, we introduce a cotriple cohomology of precrossed modules with a system of local
coefficients.

Notice that $\pcm$ is a category of interest in the sense of \cite{orz}; this follows
from the fact that $\pcm$ is equivalent to pre-cat$^1$-groups and from the tripleability
of $\pcm$ over $\set$.
\begin{lemma}\label{equiv.lem1}
    Let $\tgmu$ be a precrossed module, $A$ an abelian group. Then $(A,1,0)$ is a
    $\tgmu$-module (in the category of interest $\pcm$) if and only if $A$ is a
    $G$-equivariant $T$-module in the sense of \cite{cega} and in this case
\begin{equation*}
  \der(\tgmu,(A,1,0))\cong\der_G(T,A).
\end{equation*}
\end{lemma}
\prf If $(A,1,0)$ is a $\tgmu$-module, there is a split singular extension in $\pcm$
\begin{equation*}
  (A,1,0)\remb(T',G',\widetilde{\mu})\;\pile{\ronto\\ \leftarrow}\;\tgmu.
\end{equation*}
In particular there are split extensions of groups $A\remb T'\;\pile{\ronto\\
\leftarrow}\;T$ and $1\remb G'\;\pile{\ronto\\ \leftarrow}\;G$, so that we can assume
that $G'=G$, and that $T'\cong A\rt T$, where the action of $T$ on $A$ is by conjugation
via the splitting; we can also assume that the map $A\rt T\rw T$ is the projection and
$T\rw A\rt T$ is the inclusion. So we have the split extension in $\pcm$
\begin{equation}\label{equiv.eq1}
\begin{diagram}[scriptlabels,s=2em]
    (A,1,0)&\rEmbed & (A\rtimes T,G,\widetilde{\mu}) & \pile{\rOnto^{(\pr_T,\id_G)}\\
    \lTo_{(i,\id_G)}} & \tgmu.
\end{diagram}
\end{equation}
The action of $G$ on $A\rt T$ induces an action of $G$ on $A$; in fact, since
$(\pr_T,\id_G)$ is a map of precrossed modules, $\pr_T(\,^{g}(a,1))=1$ for all $g\in G,
\;\; a\in A$. Since the maps in the split extension (\ref{equiv.eq1}) are maps of
precrossed modules, we have, for all $a\in A,\;\;g\in G,\;\;t\in T$,
$\;\widetilde{\mu}(a,t)=\mu(t)$,
$\;^{g}(a,t)=\,^{g}(a,1)\,^{g}(0,t)=(\,^{g}a,1)(0,\,^{g}t)=(\,^{g}a,\,^{g}t)$. In
particular we obtain, for all $a\in A,\;\;g\in G, \;\;t\in T$,
$(\,^{g}(\,^{t}a),\,^{g}t)=\,^{g}(\,^{t}a,t)=\,^{g}((0,t)(a,1))=(0,\,^{g}t)(\,^{g}a,1)=
(\,^{^{g}t}(\,^{g}a),\,^{g}t)$. Hence $^{g}(\,^{t}a)=\,^{^{g}t}(\,^{g}a)$ so that, in the
terminology of \cite{cega} $A$ is a $G$-equivariant $T$-module. Conversely if $A$ is a
$G$-equivariant $T$-module , (\ref{equiv.eq1}) is a split singular extension in $\pcm$,
so $(A,1,0)$ is a $\tgmu$-module. We have
\begin{equation*}
  \der(\tgmu,(A,1,0))\cong\hund{\pcm/\tgmu}(\tgmu,(A\rtimes
  T,G,\widetilde{\mu})).
\end{equation*}
We now show that there is an isomorphism
\begin{equation*}
  \alpha:\quad\hund{\pcm/\tgmu}(\tgmu,(A\rtimes
  T,G,\widetilde{\mu}))\rw\der_G(T,A).
\end{equation*}
Let $\alpha((D,\id_T),\id_G)=D$. Notice that $D\in\der_G(T,A)$; in fact $D\in\der(T,A)$
and since $((D,\id_T),\id_G)$ is a morphism of precrossed modules,
$(D(^{g}t),^{g}t)=^{g}(D(t),t)$ so that $D(^{g}t)=^{g}D(t)$. Clearly $\alpha$ is
injective. Let $D\in\der_G(T,A)$, then $\;((D,\id_T),\id_G)\;$ is\, a\, morphism\, of\,
precrossed\, modules\, over\, $\:\tgmu$ and $\alpha((D,\id_T),\id_G)=D$ \enpr

\smallskip

Notice that if $\tgmu$ is a crossed module and $A$ is a $\p$-module, $A$ is a
$G$-equivariant trivial $T$-module in the terminology of \cite[p. 15]{cega}. If
$I:\cm\rw\pcm$ is the inclusion, we have in this case
\begin{equation*}
  \der(I\tgmu,(A,1,0))\cong\der(\tgmu,A).
\end{equation*}
\begin{remark}\label{equiv.rem1}
\end{remark}
An equivalent version of Lemma \ref{equiv.lem1} is obtained by working in
pre-cat$^1$-groups rather than in $\pcm$. Let $(A,0,0)$ and $(T\rt G,d_0,d_1)$ be the
pre-cat$^1$-groups corresponding to the precrossed modules $(A,1,0)$ and $\tgmu$
respectively. It is easily checked that $(A,0,0)$ is a $(T\rtimes G,d_0,d_1)$-module if
and only if $A$ is a $(T\rtimes G)$-module and
\begin{equation*}
  \der((T\rtimes G,d_0,d_1),(A,0,0))\cong\{D\in\der(T\rtimes G,A)\:|\;D(1,G)=0\}.
\end{equation*}
The two versions of the lemma are clearly equivalent. Recall in fact \cite[Theorem
2.2]{cega} that the categories of $G$-equivariant $T$-modules and that of $(T\rtimes
G)$-modules are equivalent, with the action of $T\rt G$ on a $G$-equivariant $T$-module
$A$ given by $^{(t,g)}a=\,^{t}(\,^{g}a)$; moreover it is straightforward to check that
there is an isomorphism
\begin{equation*}
  \alpha:\quad \der_G(T,A)\rw\{D\in\der(T\rtimes G,A)\:|\; D(1,G)=0\}
\end{equation*}
given by $\alpha(D)(t,g)=D(t),\; (t,g)\in T\rtimes G$.

\medskip

 Let $\overline{\gbb}$ be the cotriple on $\pcm$ of \cite{aria}.
\begin{proposition}\label{equiv.pro1}
    Let $\tgmu$ be a precrossed module, $A$ a $G$-equivariant $T$-module. Then for each $n\geq 0$
\begin{equation*}
  H^n\der( \overline{\gbb}_{\bullet}\tgmu,(A,1,0))\cong H^{n+1}_G(T,A).
\end{equation*}
\end{proposition}
\prf From \cite[Theorem 2.6]{cega} for each $n\geq 0$
\begin{equation*}
  H^n_G(T,A)\cong H^n(B_{T\rtimes G}, B_G,A)\cong H^n(\ker(C^*(T\rtimes
  G,A)\supar{r_{\bullet}}C^*(G,A)),
\end{equation*}
where $C^*(T\rtimes G,A)$ and $C^*(G,A)$ are the ordinary cochain complexes for computing
group cohomology and $r_{\bullet}$ are the restriction maps.

Denote $\overline{\gbb}_{\bullet}\tgmu=(T_{\bullet},G_{\bullet},\mu_{\bullet})$. From
Lemma \ref{equiv.lem1} and Remark \ref{equiv.rem1}, there is a short exact sequence for
each $n$
\begin{equation*}
  0\rw\der((T_n,G_n,\mu_n),(A,1,0))\rw\der(T_n\rtimes G_n,A)\rw \der(G_n,A)\rw 0
\end{equation*}
where the map $\der(T_n\rtimes G_n,A)\rw\der(G_n,A)$ is restriction. It is proved in
\cite[p. 12]{aria} that $T_{\bullet}\rtimes G_{\bullet}\rw T\rtimes G$ and
$G_{\bullet}\rw G$ are free simplicial resolutions. Hence, if $\bot$ is the ordinary free
cotriple on Groups we obtain short exact sequences of cochain complexes
\begin{diagram}[s=1.8em,h=1.8em]
    0 & \rTo & \der((T_{\bullet},G_{\bullet},\mu_{\bullet}),(A,1,0)) & \rTo & \der(T_{\bullet}\rtimes G_{\bullet},A) & \rTo &
    \der(G_{\bullet},A) & \rTo & 0 \\
    && \dTo && \dTo_{\wr} && \dTo_{\wr}\\
    0 & \rTo & \ker r'_{\bullet} & \rTo & \der(\bot_{\bullet}(T\rtimes G),A) & \rTo^{r'_{\bullet}} &
    \der(\bot_{\bullet}G,A) & \rTo & 0\\
    && \dTo && \dTo_{\wr} && \dTo_{\wr}\\
    0 & \rTo & \ker r_{\bullet} & \rTo & C^*(T\rtimes G,A) & \rTo^{r_{\bullet}} &
    C^*(G,A) & \rTo & 0\,.
\end{diagram}
In the above diagram the maps $\der(T_{\bullet}\rtimes
G_{\bullet},A)\!\supar{\sim}\!\der(\bot_{\bullet}(T\rtimes G),A)$,
$\;\der(G_{\bullet},A)\supar{\sim} \der(\bot_{\bullet}G,A)$ are homotopy equivalences,
and $\der(\bot_{\bullet}(T\rtimes G),A)\rw C^*(T\rtimes G,A)$, $\;
\der(\bot_{\bullet}G,A) \rw C^*(G,A)$ are the natural cochain maps of the Barr-Beck
theory which induce isomorphisms in cohomology (see \cite{bab2}). Taking the long exact
cohomology sequences in each row of the above diagram and applying the Five Lemma we
obtain for all $n>0$
\begin{equation*}
  H^n\der((T_{\bullet},G_{\bullet},\mu_{\bullet}),(A,1,0))\cong H^n\ker r'_{\bullet}\cong
  H^{n+1}\ker r_{\bullet}\cong H^{n+1}_G(T,A).
\end{equation*}
Finally, by Lemma \ref{equiv.lem1}, by \cite[p. 11]{cega} and by general properties of
cotriple cohomology
\begin{equation*}
  H^0\der((T_{\bullet},G_{\bullet},\mu_{\bullet}),(A,1,0))\cong
  \der(\tgmu,(A,1,0))\cong\der_G(T,A)\cong H^1_G(T,A).
\end{equation*}
\enpr
\begin{theorem}\label{equiv.the1}
    Let $\tgmu$ be a crossed module, $A$ a $\p$-module. Let $s_n:\overline{\gbb}^n I\Rightarrow
    I\gbb^n$, $\;I:\cm\hookrightarrow \pcm$ be as in \cite{aria} and
    $\overline{s}_n=\der(s_n,A)$. There exists a long exact cohomology sequence
\begin{equation*}
\begin{split}
   & 0\rw H^0\mathrm{coker\,}\overline{s}_{\bullet}\rw D^1(\tgmu,A)\rw H^2_G(T,A)\rw H^1
   \mathrm{coker\,}\overline{s}_{\bullet}\rw \\
   & \rw D^2(\tgmu,A)\rw H^3_G(T,A)\rw \cdots
\end{split}
\end{equation*}
\end{theorem}
\prf It is proved in \cite[p. 14]{aria} that there exists a surjective homomorphism of
resolutions
\begin{equation}\label{equiv.eq2}
\begin{diagram}[h=1.8em]
    \overline{\gbb}^n I\tgmu & \quad\cdots\quad & \overline{\gbb}^2 I\tgmu & \pile{\rTo\\ \rTo} &
    \overline{\gbb} I\tgmu & \rTo & \tgmu\\
    \dTo_{s_n} && \dTo_{s_2}&& \dTo_{s_1} && \dEq\\
    I {\gbb}^n\tgmu & \quad\cdots\quad & I{\gbb}^2\tgmu & \pile{\rTo\\ \rTo} &
    I\gbb\tgmu & \rTo & I \tgmu\;.
\end{diagram}
\end{equation}
The natural transformation $s_n:\overline{\gbb}^nI\rw I\gbb^n$ is defined inductively by
$s_{n+1}=s_n\gbb\cir\overline{\gbb}^n s_1$ for every $n\geq 1$, and $s_1$ is the natural
transformation sending each crossed module $\tgmu$ to the canonical projection
$\overline{\gbb}\tgmu\rw\mathcal{P}\overline{\gbb}\tgmu=\gbb\tgmu$ where
$\mathcal{P}:\pcm\rw\cm$ is Peiffer abelianization \cite{aria}. Hence (\ref{equiv.eq2})
gives rise to an injective morphism of cochain complexes
\begin{diagram}[h=1.8em,s=1.8em]
    0 & \rTo & \der(I \gbb\tgmu,(A,1,0)) & \rTo & \der(I \gbb^2\tgmu,(A,1,0)) & \rTo & \cdots\\
    && \dInto_{\overline{s}_1} && \dInto_{\overline{s}_2} \\
    0 & \rTo & \der(\overline{\gbb}I\tgmu,(A,1,0)) & \rTo & \der(\overline{\gbb}^2 I\tgmu,(A,1,0)) & \rTo & \cdots
\end{diagram}
Since $\der(I\gbb_{\bullet}\tgmu,(A,1,0))\cong\der(\gbb_{\bullet}\tgmu,A)$ we therefore
have a short exact sequence of cochain complexes
\begin{equation*}
  \der(\gbb_{\bullet}\tgmu,A)\overset{\overline{s}_{\bullet}}{\remb}\der(
  \overline{\gbb}_{\bullet}I\tgmu,(A,1,0))\ronto\mathrm{coker\,}\overline{s}_{\bullet}.
\end{equation*}
Since $H^0\der(\gbb_{\bullet}\tgmu,A)\cong H^0\der(\overline{\gbb}_{\bullet} I
\tgmu,A)\cong\der_G(T,A)$, taking the corresponding long exact cohomology sequence and
using Proposition \ref{equiv.pro1} the result follows. \enpr

\medskip

The following counterexample shows that in general $\;D^n(\tgmu,A)\;$ and
$H^{n+1}_G(T,A)$ are not isomorphic for $n>0$. Let $C_\infty$ be the infinite cyclic
group with generator $t$, let $H$ be the subgroup of $C_\infty$ generated by $t^2$ so
that $C_\infty/H\cong C_2$, the cyclic group of order 2. Consider the crossed module
$\Phi=(H,C_\infty,i)$ and let $A$ be a $C_2$-module. Let $\sigma$ be the generator of
$C_2$ and $N=1+\sigma$; since $H^n(C_\infty,A)=0$ for $n>1$ by Theorem \ref{ascro.the1}
\begin{equation*}
  D^n(\Phi,A)\!\cong\! H^{n+1}(C_2,C_\infty;A)\!\cong\! H^{n+2}(C_2;A)\!\cong\!
  \begin{cases}
   \!\! \{a\in A:Na=0\}/(\sigma\!-\!1)A\!\! & n\;\text{odd}, \\
   \! A^{C_2}/NA & \!\!\!\!\!\!\!\!\!\!\!\!\!\!\!\!\!\!\!\! n\; \text{even, } n>0.
  \end{cases}
\end{equation*}
On the other hand, since $H$ acts trivially on $A$, in the terminology of \cite{cega} $A$
is a $C_\infty$-equivariant trivial $H$-module, so that by \cite[Corollary 3.7]{cega}
 for all $n\geq 1$
\begin{equation*}
  H^{n+1}_{C_\infty}(H,A)\cong H^{n+1}(H\rtimes C_\infty,A).
\end{equation*}
Since the action of $C_\infty$ on $H$ is trivial, $H\rtimes C_\infty\cong H\times
C_\infty$. From the Lyndon / Hochschild -- Serre spectral sequence
\begin{equation*}
  E^{pq}_2=H^p(C_\infty,H^q(H,A))\Rightarrow H^{p+q}(H\times C_\infty,A),
\end{equation*}
since $E^{pq}_2=0$ for $p\neq 0,1$, there are exact sequences, for each $n\geq 1$
\begin{equation*}
  H^1(C_\infty,H^{n-1}(H,A))\remb H^n(H\times C_\infty,A)\ronto H^0(C_\infty,H^n(H,A)).
\end{equation*}
Since $H^n(H,A)=0$ for $n\geq 2$, $\;H^1(H,A)\cong A$, we obtain
\begin{equation*}
  H^{n+1}_{C_\infty}(H,A)\cong
  \begin{cases}
    H^1(C_\infty,A)\cong A_{C_\infty} & n=1, \\
    0 & n>1.
  \end{cases}
\end{equation*}
Hence in general $D^n(\Phi,A)\neq H^{n+1}_{C_\infty}(H,A)$ for $n>0$.

We finally remark that, since $H_i(H)=0$ for all $i\geq 2$, from \cite[p. 19]{cega} there
are isomorphisms $H^{n+1}_{C_\infty}(H,A)\cong\ext^n_{C_\infty}(H_{ab},A)$ for all $n\geq
0$. Hence this counterexample also shows that despite the isomorphism
$D^0(\tgmu,A)\cong\hund{G}(T,A)\cong\hund{G\mi\mathrm{Mod}}(T_{ab},A)$ for any crossed
module $\tgmu$, in general the cohomology groups $D^n(\tgmu,A)$ and
$H^n_{G\mi\mathrm{Mod}}(T_{ab},A)$ are not isomorphic for $n>0$.

\renewcommand{\refname}{References}

\end{document}